\pdfoutput=1
\pdfoutput=1
\documentclass[final,1p,times]{elsarticle}

\usepackage{epsfig}
\usepackage{float}

\usepackage{amssymb}

\usepackage{amsmath, amssymb,amsfonts}
\usepackage[nodots]{numcompress}

\usepackage{lineno}

\newcommand{\pt}{\,\partial_t\,}
\newcommand{\px}{\,\partial_x\,}

\newcommand{\pz}{\,\partial_z\,}

\newcommand{\F}{\,\mathcal F\,}

\renewcommand{\P}{\,\mathcal P\,}
\renewcommand{\H}{\,\mathcal H\,}
\renewcommand{\S}{\,\mathcal S\,}
\newcommand{\V}{\,\mathcal V\,}
\newcommand{\dt}{\Delta t }
\newcommand{\dx}{\Delta x }

\newcommand{\dz}{\Delta z}

\begin{document}

\begin{frontmatter}



\title{A study of a WENO-TVD finite volume scheme for the numerical simulation of atmospheric advective and convective phenomena}

\author[dm]{Dante Kalise}
\address[dm]{Dipartimento di Matematica, Universit\`a di Roma ``La Sapienza'', P. Aldo Moro 2, 00185 Roma,
Italy.}
\ead{kalise@mat.uniroma1.it}

\begin{abstract}
We present a WENO-TVD scheme for the simulation of atmospheric phenomena. The scheme considers a spatial discretization via a second-order TVD flux based upon a flux-centered limiter approach, which makes use of high-order accurate extrapolated values arising from a WENO reconstruction procedure. Time discretization is performed with a third order RK-TVD scheme, and splitting is used for the inclusion of source terms. We present a comprehensive performance study of the method in atmospheric applications involving advective and convective motion. We present a set of tests for space-dependent linear advection, where we assess convergence and robustness with respect to the parameters of the scheme. We apply the method to approximate the 2D Euler equations in a series of tests for atmospheric convection.
\end{abstract}


\end{frontmatter}


\section{Introduction}
Advection and convection, understood as the class of phenomena related to fluid flow motion or transport is a concept covering different relevant situations in atmospheric modelling. The convention in the literature is to use the term \textsl{advection} when a quantity experiences motion due to the presence of an acting velocity field, which is most often related to horizontal motion, while \textsl{convection} refers to motion caused by thermodynamic considerations, primarily occurring in the vertical direction. In this article we are concerned with the development of an accurate scheme able to handle both types of motion. Models which describe such behavior are, in a first step, linear scalar advection models, including space-dependent velocity fields and, in a more elaborated formulation, the set of Euler equations for gas dynamics, parameterized in scales that are typical of atmospheric phenomena.
Achieving an accurate and physically meaningful numerical approximation of such models is undoubtedly, a challenging task. During the last decades, a method which has gained popularity among the atmospheric modelling community is the finite volume (FV) framework \cite{torolibro,levequelibro} . In a computationally efficient way, it preserves many aspects of the underlying physics (as conservation for instance), while allowing formulations that are able to achieve a high level of accuracy, which is fundamental for numerical weather prediction.
In the FV context there is a considerable amount of available methods oriented to the solution of hyperbolic system of conservation laws. We are concerned with a particular class of those methods, the so-called WENO (weighted essentially non-oscillatory methods), which is a class of methods allowing the generation of high-resolution approximations in space via a polynomial reconstruction procedure from cell averages \cite{liu,shu1}; this technique is combined with suitable fluxes and time marching procedures in order to generate a method of global high-order of accuracy. Moreover, at every step, by adequately enforcing the concept of total variation diminishing (TVD), the method avoids the production of spurious oscillations and preserves monotonicity.
This article addresses every step in the construction of the above described scheme. Once the WENO reconstruction of the variables is performed, high-order extrapolated values are available for calculation of numerical fluxes across the cell interfaces. Based upon \cite{titarevweno}, we make use of a second-order TVD flux based on a centered limiter (FLIC) approach \cite{torobillett}; we highlight the use of the centered approach instead of classical upwinding considerations \cite{godunov}. After the spatial discretization routine is completed, the scheme moves forward in time via a Runge-Kutta TVD  method \cite{shutvd} with ensures stability, preservation of the spatial accuracy and avoids the generation of new extremal values.
An important part of this work is devoted to the numerical validation of the proposed scheme. We first address 2D models of space-dependent advection where several features of the scheme are tested, to then study convective phenomena based upon the Euler equations, with a set of well-known tests for atmospheric modelling. The performance of the scheme is assessed in terms of accuracy, its ability to preserve monotonicity in the presence of sharp solutions, its robustness with respect to flux parameters, correct front locations and energy conservation.

The article is  structured as follows. In section 2 we present the full numerical WENO-TVD scheme. In section 3 we develop a comprehensive study of the method for advective models, while section 4 is devoted to the analysis of convective phenomena. Final remarks are discussed in section 5.

\section{A numerical scheme for the system of balance laws}
In this section we present a finite volume scheme for  a two-dimensional system of balance laws of the form
\begin{equation}
\pt Q+\px \F(Q)+ \pz \H(Q)=\S(Q),
\end{equation}
where $Q$ is a vector of conserved variables, $\F$ and $\H$ are physical fluxes, and $\S$ is a source term. We first indicate that our strategy will be based in a splitting scheme, as it is suggested in \cite{TOROCHEN} given the flux choice that we will make. Thus, we will first establish a numerical scheme for the system of conservation laws
\begin{equation}\label{scl}
\pt Q+\px \F(Q)+ \pz \H(Q)=0,
\end{equation}
to be combined with a procedure for the resolution of the source term dynamics
\begin{equation}\label{st}
\pt (Q)=\S(Q).
\end{equation}
In order to approximate eq. (\ref{scl}) we begin by meshing the spatial domain $\Omega_{x,z}$ into uniform control volumes $\Omega_{i,j}=[x_{i-1/2},\,x_{i+1/2}]\times [z_{j-1/2},\,z_{j+1/2}]$ of size $\dx\,\dz$; inside every control volume we average with respect to $x$ and $z$ leading to the semi-discrete scheme
\begin{equation}\label{sclh}
\frac{d Q_{i,j}(t)}{dt}=-\frac{1}{\dx}(F_{i+1/2,j}-F_{i-1/2,j})-\frac{1}{\dz}(H_{i,j+1/2}-H_{i,j-1/2})\equiv L_{i,j}(Q),
\end{equation}
where
\begin{equation}
Q_{i,j}=\frac{1}{\dx}\frac{1}{\dz}\int_{x_{i-1/2}}^{x_{i+1/2}}\,\int_{z_{j-1/2}}^{z_{j+1/2}}Q(x,z,t)\,dz\,dx,
\end{equation}
\begin{equation}\label{fluxes}
F_{i+1/2,j}=\frac{1}{\dz}\int_{z_{j-1/2}}^{z_{j+1/2}}\F(Q(x_{i+1/2},z,t))\,dz\,,\quad
H_{i,j+1/2}=\frac{1}{\dx}\int_{x_{i-1/2}}^{x_{i+1/2}}\H(Q(x,z_{j+1/2},t))\,dx.
\end{equation}
We approximate the expressions in (\ref{fluxes}) by conventional Gaussian quadrature formulas
\begin{equation}\label{fluxesh}
F_{i+1/2,j}\approx\frac{1}{2}\sum_{z_{Gp}}w_{z_{Gp}}\F(Q(x_{i+1/2},z_{Gp},t))\,dz\,,\quad
H_{i,j+1/2}\approx\frac{1}{2}\sum_{x_{Gp}}w_{x_{Gp}}\H(Q(x_{Gp},z_{j+1/2},t))\,dz,
\end{equation}
where $x_{Gp}$ and $z_{Gp}$ are prescribed Gauss points with corresponding weights $w_{x_{Gp}}$ and $w_{z_{Gp}}$ respectively.
The computation of (\ref{fluxesh}) is performed via a high-resolution approach that makes use of a WENO reconstruction procedure; after this step is completed, a polynomial of prescribed order is obtained at every cell, and therefore, at every cell interface, accurate flux calculations can be performed by taking extrapolated boundary values.

We briefly describe the WENO reconstruction procedure that is used in this article; we opted for the technique described in \cite{balsara1,balsara2} in its third order (quadratic reconstruction) version. This technique makes extensive use of the structure of the reconstruction procedure in one dimension, adding some additional mixed terms (``cross terms'') that are efficiently computed by reduced stencils. It is an optimal and easy way to implement the algorithm for achieving high-order reconstructions in 2 and 3 dimensions; it also defines an unique polynomial in every cell, which is particularly useful when space dependent source terms such as viscosity are considered.
At a given time $t$ (the subscript indicating time is omitted throughout this derivation), given the set of averaged values $\{Q_{i,j}\}$ for the whole domain, at every cell, the reconstruction procedure seeks a quadratic expansion upon a linear combination of Legendre polynomials rescaled in local coordinates $(x,z)=[-1/2,\,1/2]\times [-1/2,\,1/2]$ expressed in the form
\begin{equation}\label{rec}
Q(x,z)=Q_0+Q_xP_1(x)+Q_{xx}P_2(x)+Q_zP_1(z)+Q_{zz}P_2(z)+Q_{xz}P_1(x)P_1(z),
\end{equation}
\begin{equation}
P_1(x)=x \qquad P_2(x)=x^2-\frac{1}{12}.
\end{equation}
Except for the last term in (\ref{rec}), every coefficient can be computed by performing a dimension-by-dimension reconstruction, which we now illustrate. We assign the subscript ''0'' to the cell where we are computing the coefficients, other values indicating location and direction with respect to $Q_0$ (note that the notation is coherent with the fact that the first coefficient in the expansion $Q_0$, holds $Q_0=Q_{ij}$, i.e., the centered value). Next, for this particular problem we define three stencils
\begin{equation}
S^1=\{Q_{-2},Q_{-1},Q_{0}\}\,,\quad S^2=\{Q_{-1},Q_{0},Q_{1}\}\,,\quad S^3=\{Q_{0},Q_{1},Q_{2}\}\,,
\end{equation}
and in every stencil we compute a polynomial of the form
\begin{equation}
Q^{(i)}(x)=Q_0^{(i)}+Q_x^{(i)}P_1(x)+Q_{xx}^{(i)}P_2(x)\qquad i=1,2,3.
\end{equation}
The coefficients are given by
\begin{eqnarray}
S^{1}&:&Q^{(1)}_x=-2Q_{-1}+Q_{-2}/2+3Q_0/2,\quad Q^{(1)}_{xx}=(Q_{-2}-2Q_{-1}+Q_0)/2\,,\\
S^{2}&:&Q^{(2)}_x=(Q_1-Q_{-1})/2,\quad Q^{(2)}_{xx}=(Q_{-1}-2Q_{0}+Q_1)/2\,,\\
S^{3}&:&Q^{(3)}_x=-3Q_{0}/2+2Q_{1}-Q_{2}/2,\quad Q^{(3)}_{xx}=(Q_{0}-2Q_{-1}+Q_2)/2.
\end{eqnarray}
For every polynomial we calculate a smoothness indicator defined as
\begin{equation}
IS^{(i)}=\left(Q^{(i)}_x\right)^2+\frac{13}{3}\left(Q^{(i)}_{xx}\right)^2\,,
\end{equation}
leading to the following WENO weights:
\begin{equation}
\omega^{(i)}=\frac{\alpha^{(i)}}{\sum_{i=1}^3\alpha^{(i)}}\,,\quad \alpha^{(i)}=\frac{\lambda^{(i)}}{(\epsilon+IS^{(i)})^r}\,^,
\end{equation}
where $\epsilon$ is a parameter introduced in order to avoid division by zero; usually $\epsilon=10^{-12}$. The scheme is rather insensitive to the parameter $r$, which we set $r=5$. The parameter $\lambda$ is usually computed in an optimal way to increase the accuracy of the reconstruction at certain points; we opt for a centered approach instead, thus $\lambda^{(1)}=\lambda^{(3)}=1$, while $\lambda^{(2)}=100$. The 1D reconstructed polynomial is given by
\begin{equation}
Q(x)=\omega^{(1)}Q^{(1)}(x)+\omega^{(2)}Q^{(2)}(x)+\omega^{(3)}Q^{(3)}(x).
\end{equation}
Next, a 1D reconstruction in the $z$ direction is performed in a totally analogous way. Finally, we address the computation of the mixed term $Q_{xz}$, which is calculated in a 2D fashion. Keeping the same convention regarding location subscripts as in 1D, \cite{balsara1} considers 4 formulas for the cross term upon taking all the moments around the cell. The expressions for the cross term are:
\begin{eqnarray}
Q_{xz}^{(1)}&=&Q_{1,1}-Q_{0,0}-Q_x-Q_z-Q_{xx}-Q_{zz},\\
Q_{xz}^{(2)}&=&-Q_{1,-1}+Q_{0,0}+Q_x-Q_z+Q_{xx}+Q_{zz},\\
Q_{xz}^{(3)}&=&-Q_{-1,1}+Q_{0,0}-Q_x+Q_z+Q_{xx}+Q_{zz},\\
Q_{xz}^{(4)}&=&Q_{-1,-1}-Q_{0,0}+Q_x+Q_z-Q_{xx}-Q_{zz},
\end{eqnarray}
and the corresponding smoothness indicators are given by
\begin{equation}
IS^{(i)}=4\left(Q^{(i)}_{xx}\right)^2+4\left(Q^{(i)}_{zz}\right)^2+\left(Q^{(i)}_{xz}\right)^2.
\end{equation}
Note that in the first part of the reconstruction, when the weights were computed, a larger suboptimal weight was assigned to the central stencil, which is a way to ensure stability and robustness of the algorithm by sacrificing additional order in the approximation (for more details, see \cite{dumbserkaser}). However, for this term, the numerators assigned to the corresponding $\alpha$'s remains the same for every expression. The computation of this term concludes the reconstruction procedure, and now we have at our disposal one polynomial per cell that can be used to calculate values at the boundaries or inside the cell. The next step in our numerical scheme consists of the calculation of the numerical fluxes (\ref{fluxesh}), which will use extrapolated boundary values of the reconstructed polynomials. Rather than the use of the classical WENO scheme (as in \cite{liu} or \cite{shu1}), which performs this calculation via a first order flux, we opt for the WENO-TVD approach described in \cite{titarevweno}. We make use of the 2D extension of the flux-limiter-centred scheme (FLIC) approach presented in \cite{torobillett,torolibro}, which is a second-order, centered and non-oscillatory flux. In our case, it consists of a flux-limited version of a generalized Lax Wendroff flux, using  as a low-order flux the GFORCE (generalized first order centred) flux \cite{toroforce}, which can be interpreted as a convex combination of Lax-Friedrichs and Lax-Wendroff-type of fluxes:
\begin{equation}
F_{i+1/2,j}^{FLIC}=F_{i+1/2,j}^{GFORCE}+\psi_{i+1/2,j}\left(F_{i+1/2,j}^{LW}-F_{i+1/2,j}^{GFORCE}\right),
\end{equation}

where
\begin{eqnarray}
F_{i+1/2,j}^{GFORCE}&=&F_{i+1/2,j}^{GFORCE}\left(Q_{i+1/2,j}^L,Q_{i+1/2,j}^R\right)=\omega F_{i+1/2,j}^{LW}+(1-\omega)F_{i+1/2,j}^{LF}\,,\\
F_{i+1/2,j}^{LF}&=&\frac12\left(\F\left(Q_{i+1/2,j}^L\right)+\F\left(Q_{i+1/2,j}^R\right)-\frac12\frac{\dx}{\dt}\left(Q_{i+1/2,j}^R-Q_{i+1/2,j}^L\right)\right)\,,\\
F_{i+1/2,j}^{LW}&=&\F\left( Q_{i+1/2,j}^* \right)\,,\\
Q_{i+1/2,j}^*&=&\frac12\left(Q_{i+1/2,j}^L+Q_{i+1/2,j}^R\right)-\frac{\dt}{\dx}\left(\F\left(Q_{i+1/2,j}^R\right)-\F\left(Q_{i+1/2,j}^L\right)\right).
\end{eqnarray}
The parameter $\omega$ varies between 0 and 1, and is chosen in a compatible manner with the CFL number in order to ensure monotonicity.
We have omitted the formulas for the remaining cell boundaries, but they can be derived in a straightforward manner. Also note that even though the formulas are written along the boundary '$i+1/2,j$', the use of the Gaussian quadrature formula will replace the axes'$j$' by Gauss points and therefore this subscript must be understood in that sense.
It is important to notice that so far we are deriving expressions for the semi-discrete approximation of the system of conservation laws, however, the fluxes include the parameter $\dt$ which arises from the averaging operators that originate these fluxes. Thus, in the spatial discretization of the system, the time stepping enters just as a parameter. At the end of the derivation of the scheme, when we present the time discretization of eq. (\ref{sclh}), $\dt$ will be considered as ``marching parameter'' in the sense that its inclusion in the formulas will generate an updated state in time.

The function $\psi_{i+1/2,j}=\psi_{i+1/2,j}(r^L_{i+1/2,j},r^R_{i+1/2,j})$ is a flux limiter; a slight variation of the usual limiters has to be considered in this context since we use a centered flux instead of an upwind approach (the reader can refer to \cite[Ch 13.]{torolibro} for more details); in our case we mainly use the SUPERBEE limiter, which on its centered version reads:
\begin{equation}
\psi(r)=\begin{cases}
0 &\mbox{if } r\leq 0,\\
2r &\mbox{if } 0\leq r \leq \frac12,\\
1 &\mbox{if } \frac12\leq r \leq 1,\\
\min \left\{2,\phi_g+(1-\phi_g)r\right\} r\geq 1,
\end{cases}
\qquad \phi_g=\frac{1-|c|}{1+|c|},
\end{equation}
where $c$ corresponds to the Courant number which depends on the problem. The limiter depends on the flow parameter $r$, which will be defined upon a physical quantity $e$ of the system. Once $e$ has been obtained from the discretized variables, left and right flow parameters are given by
\begin{equation}\label{limiter}
r^L_{i+1/2,j}=\frac{e^R_{i-1/2,j}-e^L_{i-1/2,j}}{e^R_{i+1/2,j}-e^L_{i+1/2,j}}\,,\qquad r^R_{i+1/2,j}=\frac{e^R_{i+3/2,j}-e^L_{i+3/2,j}}{e^R_{i+1/2,j}-e^L_{i+1/2,j}},
\end{equation}
and finally,
\begin{equation}
\psi_{i+1/2,j}=\min(\psi(r^L_{i+1/2,j}),\psi(r^R_{i+1/2,j})).
\end{equation}

The above described procedure starts with a set of averaged values and ends with a numerical approximation of the space operators involved in eq. (\ref{scl}). The resulting scheme is still continuous in time, and we conclude this section by discretizing this operator in a manner that is consistent with the choices that we have made in the generation of the space discretization operator. At a given starting time $t^n$, we begin by considering the semi-discrete scheme
\begin{equation}
\frac{d Q_{i,j}(t)}{dt}= L_{i,j}(Q),
\end{equation}
bringing the system to a final state $t^{n+1}$ with a time stepping $\dt$. In order to preserve high-order and non-oscillatory properties in time, we consider the well-known family of explicit TVD Runge-Kutta schemes \cite{shutvd}, in particular its third order version
\begin{eqnarray}
Q_{i,j}^{n+\frac13}&=&Q_{i,j}^n+\dt \,L_{i,j}(Q_{i,j}^n),\\
Q_{i,j}^{n+\frac23}&=&\frac34Q^n_{i,j}+\frac14Q_{i,j}^{n+\frac13}+\frac14\dt\, L_{i,j}(Q_{i,j}^{n+\frac13}), \\
Q_{i,j}^{n+1}&=&\frac13 Q^n_{i,j}+\frac23 Q_{i,j}^{n+\frac23} +\frac23\dt\, L_{i,j}(Q_{i,j}^{n+\frac23}). \\
\end{eqnarray}

We end this section with the inclusion of the source term. The source term appearing in eq. (\ref{st}), in the simplest case will not depend on space nor space derivatives, and therefore it can be averaged in space and solved in the same manner as the above presented time discretization, by replacing $L_{i,j}(Q_{i,j}^n)$ by $\S(Q_{i,j}^n)$. If we denote by the $\mathfrak{L}(\dt)$ the fully discrete operator that brings the system of conservation laws (\ref{scl}) $\dt$ units ahead in time, and by $\mathfrak{S}(\dt)$ the fully discrete operator that updates the source term (\ref{st}) in $\dt$ units, we preserve, at least, second order accuracy in time by implementing a Strang splitting \cite{strang} in the form
\begin{equation}
Q^{n+1}_{i,j}=\mathfrak{S}(\dt/2)\mathfrak{L}(\dt)\mathfrak{S}(\dt/2)Q_{i,j}^n.
\end{equation}
If the source term does depends either on space or space derivatives, such a viscosity for instance, the averaging procedure will require the evaluation of the source term integral
\begin{equation}
S_{i,j}=\frac{1}{\dx}\frac{1}{\dz}\int_{x_{i-1/2}}^{x_{i+1/2}}\,\int_{z_{j-1/2}}^{z_{j+1/2}}\S(Q(x,z,t))\,dz\,dx.
\end{equation}
Proceeding in the same way as we did for the fluxes, we approximate this integral by a suitable double Gaussian quadrature,
\begin{equation}
S_{i,j}\approx\frac{1}{4}\sum_{x_{Gp}}\sum_{z_{Gp}}w_{x_{Gp}}w_{z_{Gp}}\S(Q(x_{Gp},z_{Gp},t)),
\end{equation}
where we make use of the same reconstruction procedure previously described in order to obtain values of $Q$ inside every cell.

\section{Advective tests: Linear advection with space-dependent coefficients}
In this section we implement three different test cases based on the equation
\begin{equation}\label{adveccion}
\pt Q + \px(a(x,z)Q)+\pz(b(x,z)Q)=0,
\end{equation}
which describes the evolution of a scalar quantity $Q$ through a 2D domain $\Omega$, holding suitable initial and boundary conditions. The tests for this equation aims to recover the theoretically expected second-order for the convergence rate of the scheme, to study variations on the weighting parameter $\omega$ in the flux and effects of the limiter choice for sharp initial conditions. Test settings varies from one test to the other; though, one common aspect is the computation of the time stepping $\dt$, which once mesh parameters, $\dx$ and $\dz$, and the background flow $(a,b)$ have been specified, is computed via
\begin{equation}
\dt=CFL\,\,\min\left(\frac{\dx}{\max_{\Omega}\,|a|}, \frac{\dz}{\max_{\Omega}\,|b|}\right),
\end{equation}
with $CFL$ the classical Courant number which by default is set to $CFL=0.45$, consistently chosen for a value of $\omega=0.5$ (there is a direct relation between the weighting parameter and the CFL number in order to preserve monotonicity of the scheme; see \cite{toroforce} for precise details). A last common aspect the choice of flow parameter $r$, which is computed by taking $e=Q\,$ in eqs. (\ref{limiter}).
\subsection{2D Linear advection with constant coefficients}
The first case that we consider is a 2D linear advection equation with constant coefficients. We set $a=1$, $b=1$, $\Omega=[0,\,1]^2$, and an initial profile given by
\begin{equation}
Q(x,z,0)=\sin(2\pi x)\sin(2\pi z),
\end{equation}
together with periodic boundary conditions. We advect the initial profile for 10 periods, i.e., final time of the simulation is $t=10$ [s]. Results for the convergence rates are shown in table~\ref{tab:t1}; expected second-order for smooth solutions is achieved in $L_1$ and $L_{\infty}$ discrete norms, while an additional order is obtained for the $L_1$ norm, which is given by the fact that we are using a third order reconstruction in space.

\begin{table}[H]
\centering
\caption{Convergence rates for the 2D linear advection problem with constant coefficients after $t=10[s]$ (10 periods).}\label{tab:t1}\vskip 3mm
\begin{tabular}{c c c c c}
\hline\\
N & $L_{\infty}$ error &  $L_{\infty}$ order & $L_1$ error &  $L_1$ order \\ [0.5ex]
\hline\\
50 &  1.2637e-002 &  & 2.005e-002 & \\
100 &   2.4316e-003 & 2.4 & 2.5911e-003 & 3.0 \\
200 &  6.1038e-004 & 2.0 & 3.5190e-004 & 2.9\\
400 &  1.4790e-004 & 2.0  & 4.9012e-005 & 2.9 \\
\hline
\end{tabular}
\end{table}

\subsection{Swirling flow}
\begin{figure}
\caption{Swirling flow test problem. Initial condition (left) and maximal flow deformation at $t=2.5$ [s] with $200\times 200$ elements (right).}
\centering
\includegraphics[width=\textwidth]{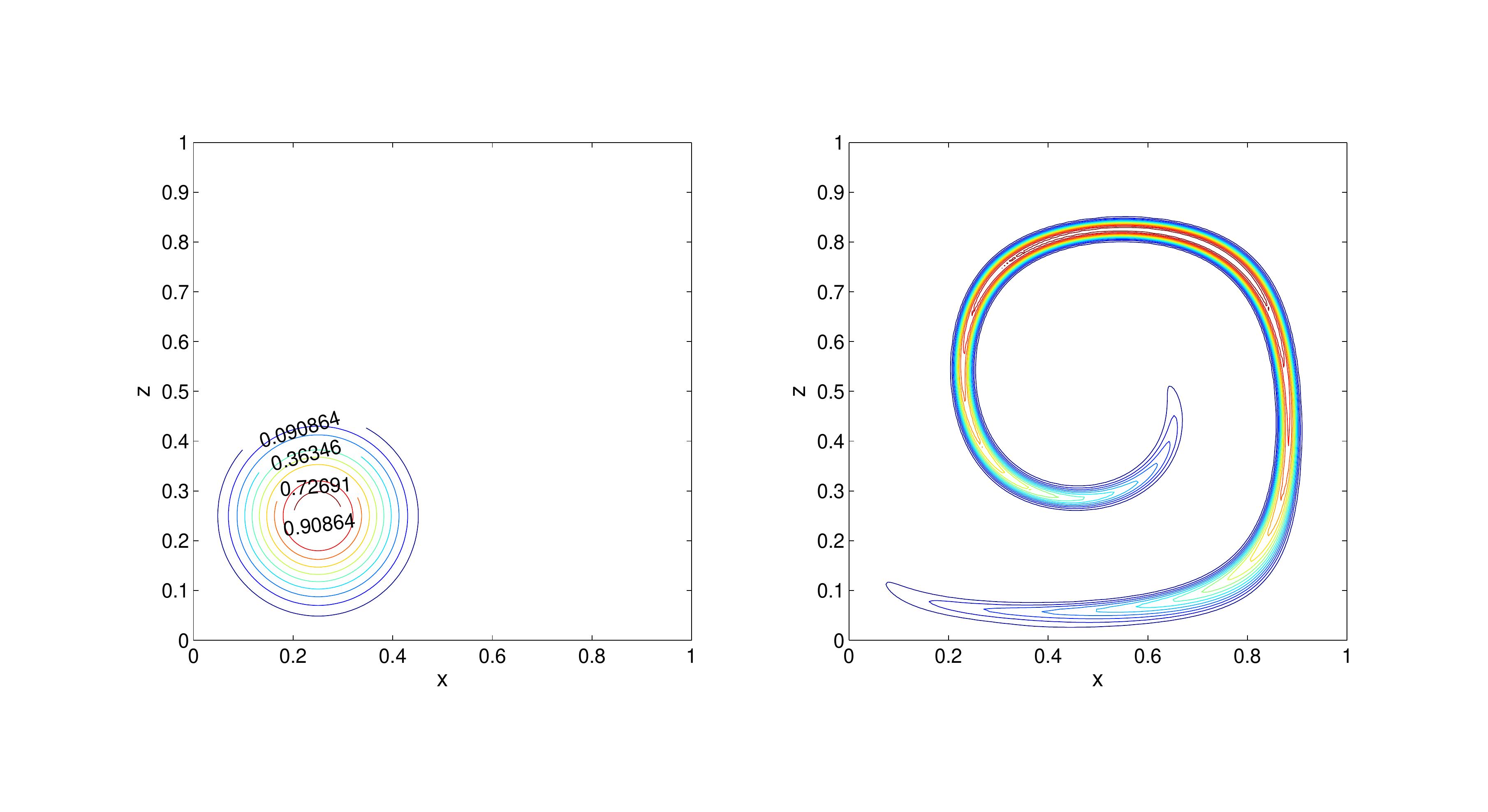}
\label{swirlinginit}
\end{figure}

The second test that we present is a more stringent case, proposed in \cite{leveque}. We make use of the setting proposed in \cite{durran}: the velocity field is given by
\begin{equation}
a=\sin^2(\pi x)\sin(2\pi z)g(t),\qquad b=-\sin^2(\pi z)\sin(2\pi x)g(t),
\end{equation}
\begin{equation}
g(t)=cos(\pi t/T).
\end{equation}
Again, $\Omega$ is set to be the unit square, and the initial condition is taken as
\begin{equation}
Q(x,z,0)=\frac12(1+\cos (\pi r),\qquad r=\min (1,4\sqrt{(x-0.25)^2+(z-0.25)^2}).
\end{equation}
Note that the velocity field vanishes at the boundary of the domain. Setting final time $T=5$, we expect a maximal flow deformation at $T/2$ while, ideally, $Q(x,z,T)=Q(x,z,0)$. Figure \ref{swirlinginit} shows both exact initial (final) profile and maximal flow deformation at $T/2$ with $\dx=\dz=0.005$; the scheme preserves positivity of the initial profile at any time during simulation. Again, table \ref{tab:t2} shows that second order of accuracy is reached as expected, although the convergence requires a larger number of elements; this is most likely due to the highly deformational nature of the flow.

\begin{table}[H]
\centering
\caption{Convergence rates for the swirling flow problem at  $t=5$ [s].}\label{tab:t2}\vskip 3mm
\begin{tabular}{c c c c c}
\hline\\
N & $L_{\infty}$ error &  $L_{\infty}$ order & $L_1$ error &  $L_1$ order \\ [0.5ex]
\hline\\
50 &  7.1093e-001 &  & 1.0719e-000 & \\
100 &   5.2363e-001 & 0.4 & 6.5837e-001 & 0.7 \\
200 &  2.4912e-001 & 1.1 & 2.2685e-001 & 1.5\\
400 &  4.5618e-002 & 2.4  & 3.5510e-002 & 2.7 \\
\hline
\end{tabular}
\end{table}

\begin{figure}
\caption{Swirling flow test problem. Low-resolution experiments with $100\times 100$ elements at $t=5$ [s]. Varying values of $\omega$ in the GFORCE flux: From left to right, from top to bottom: $\omega=0.25,\,0.5,\,0.75,\,0.9$.}
\begin{minipage}[b]{0.5\linewidth}
\centering
\includegraphics[scale=0.3]{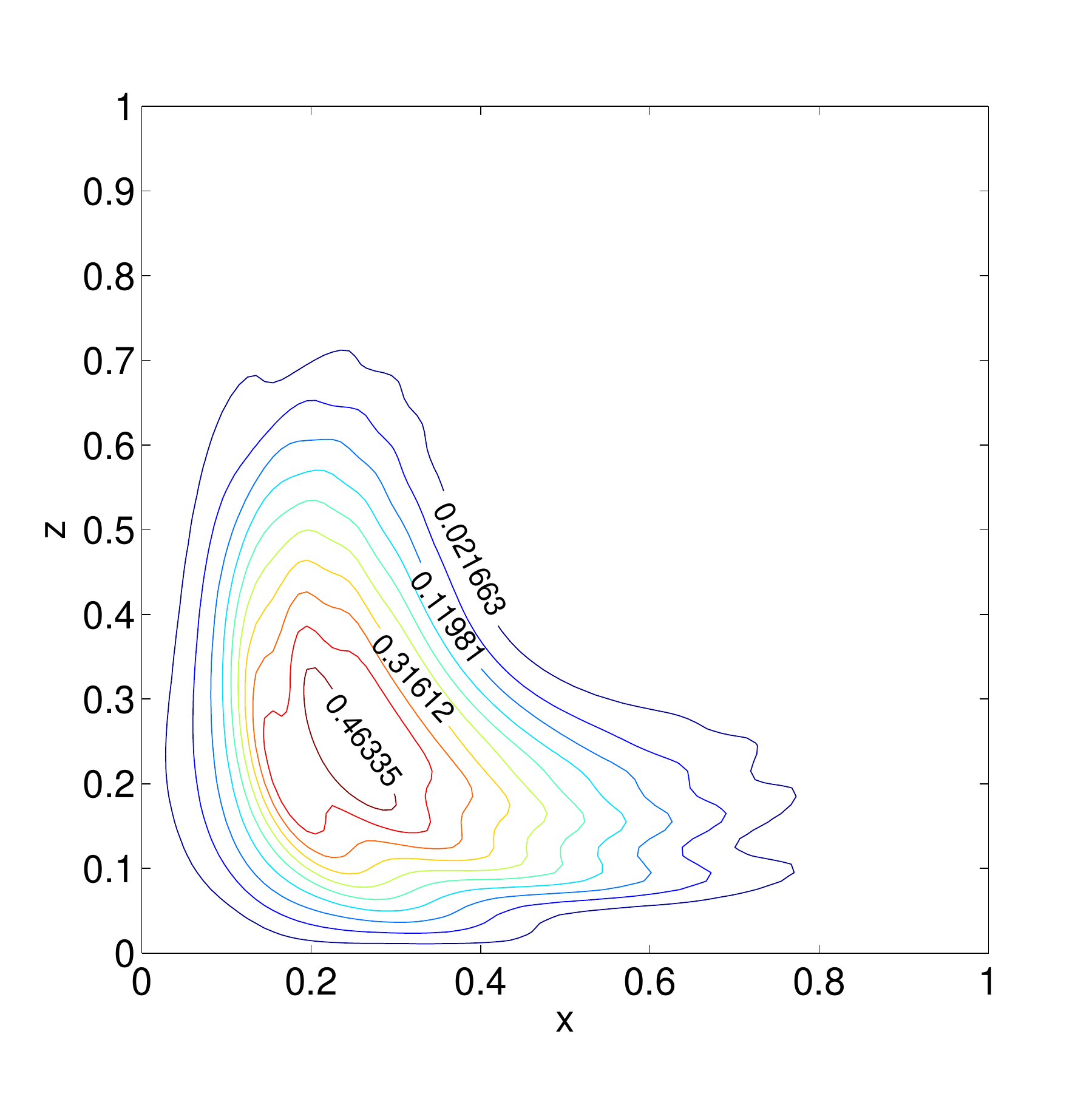}
\includegraphics[scale=0.3]{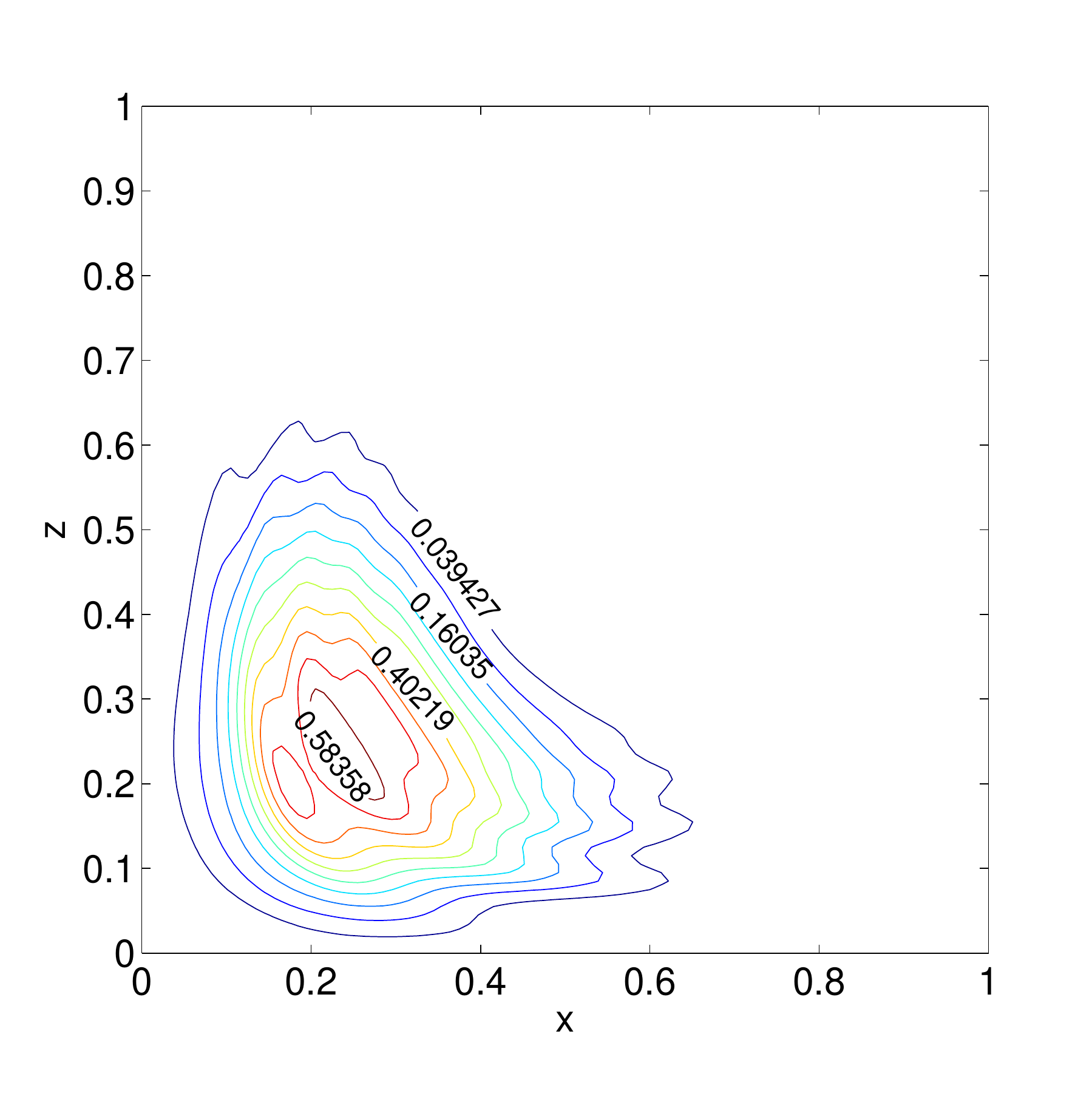}
\end{minipage}
\hspace{0.0cm}
\begin{minipage}[b]{0.5\linewidth}
\centering
\includegraphics[scale=0.3]{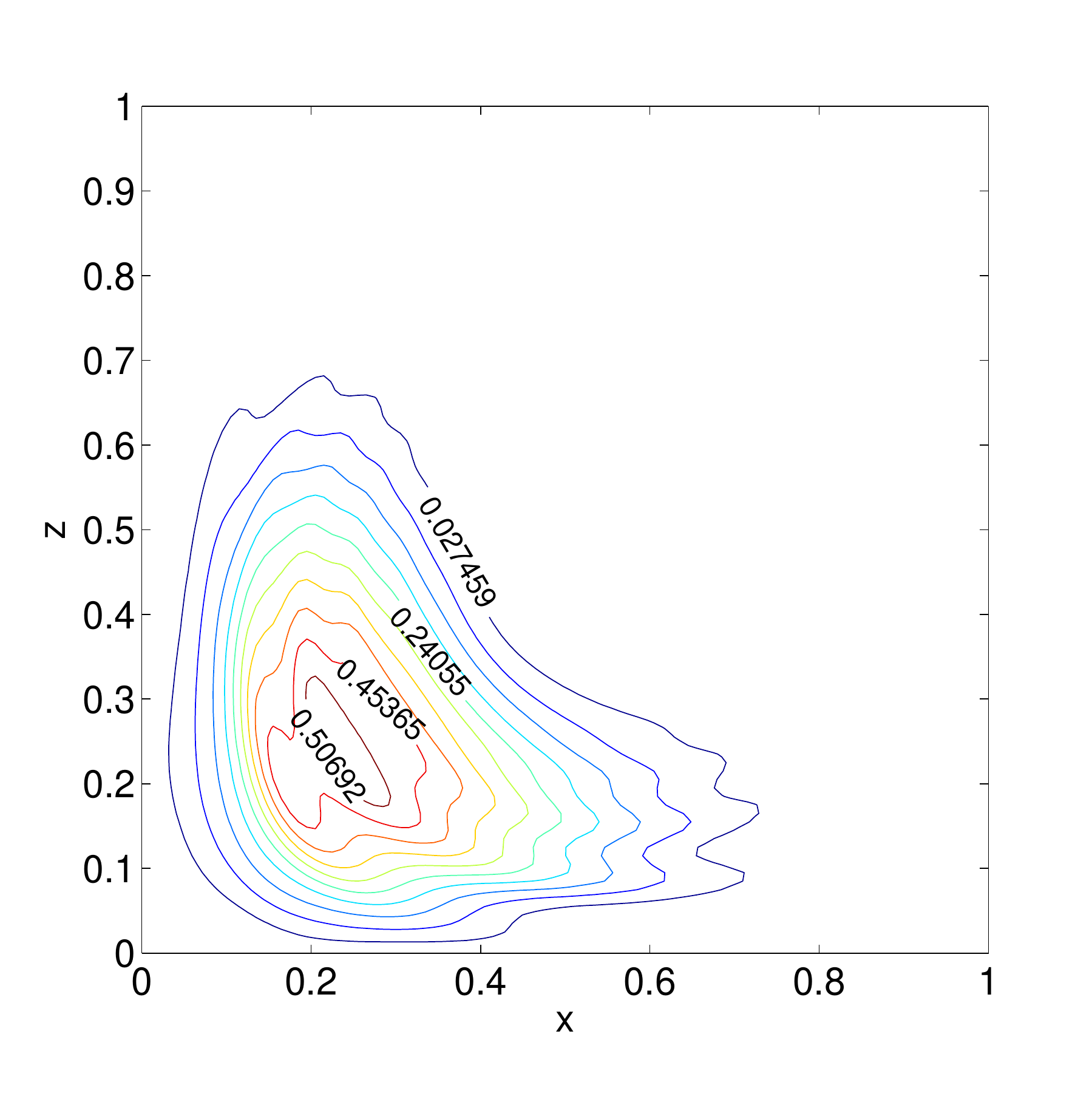}
\includegraphics[scale=0.3]{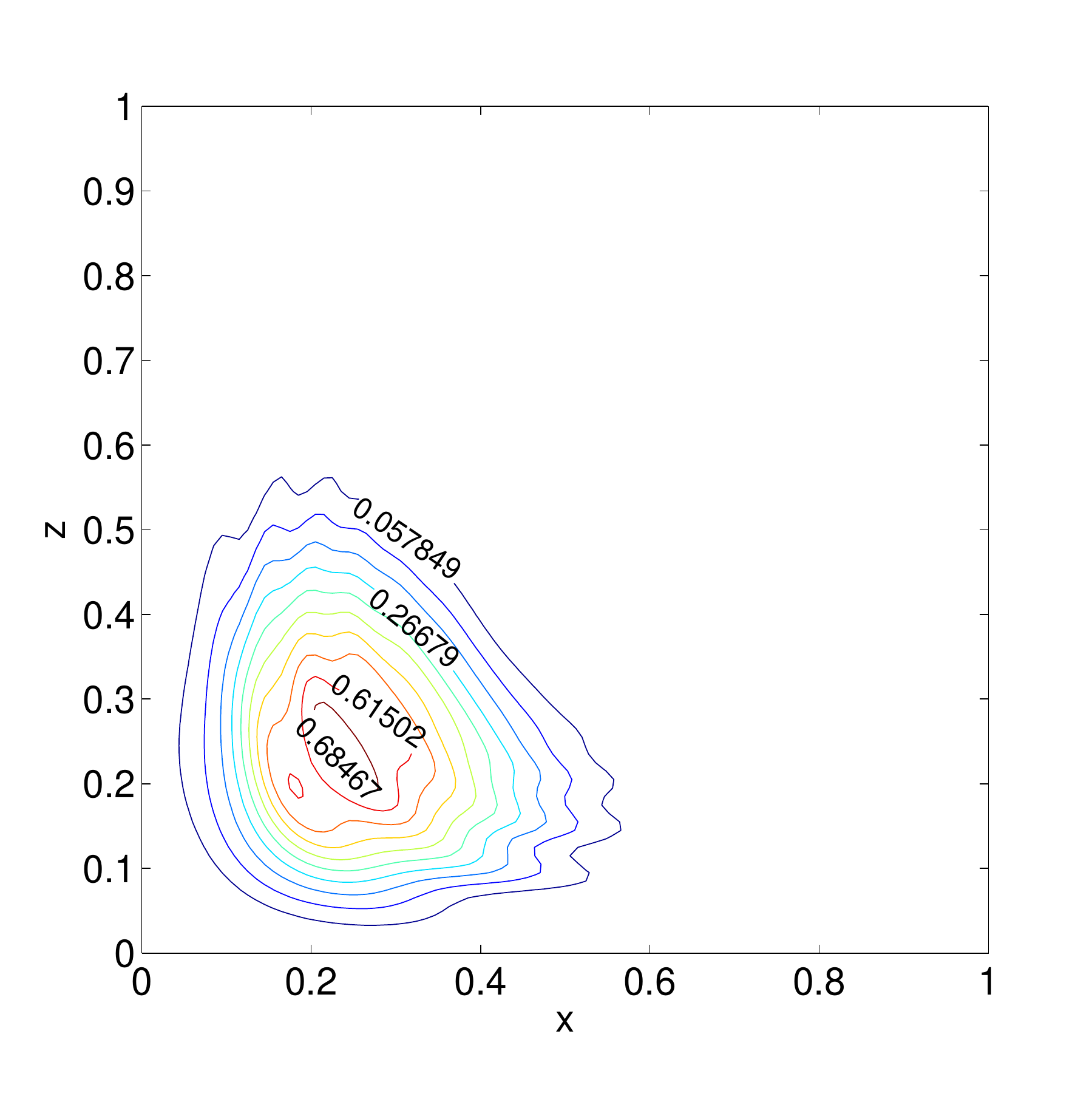}
\end{minipage}
\label{swirlinglow}
\end{figure}

\begin{figure}
\caption{Swirling flow test problem. High-resolution experiments with $200\times 200$ elements at $t=5$ [s]. Varying values of $\omega$ in the GFORCE flux: From left to right, from top to bottom: $\omega=0.25,\,0.5,\,0.75,\,0.9$.}
\begin{minipage}[b]{0.5\linewidth}
\centering
\includegraphics[scale=0.3]{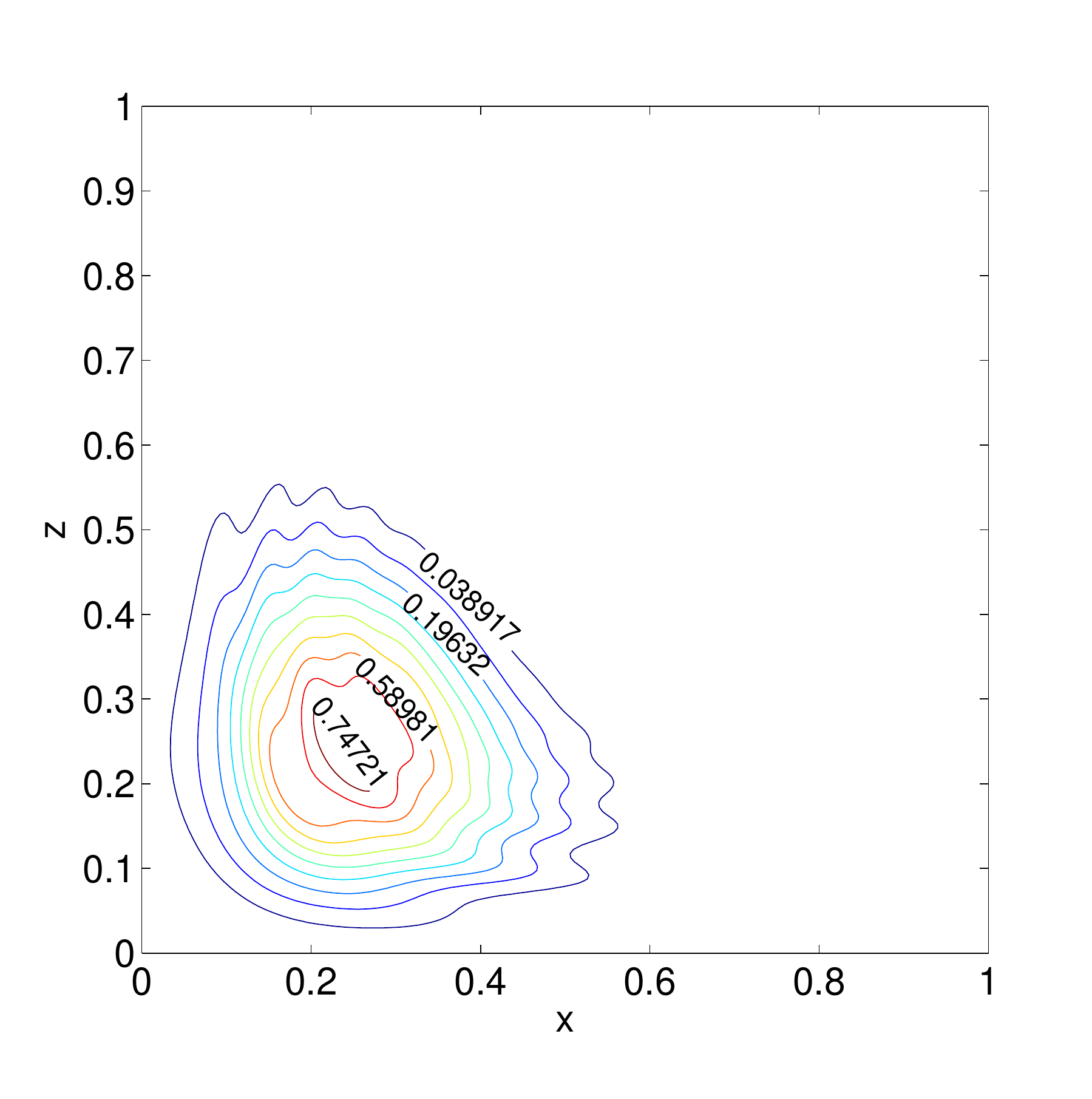}
\includegraphics[scale=0.3]{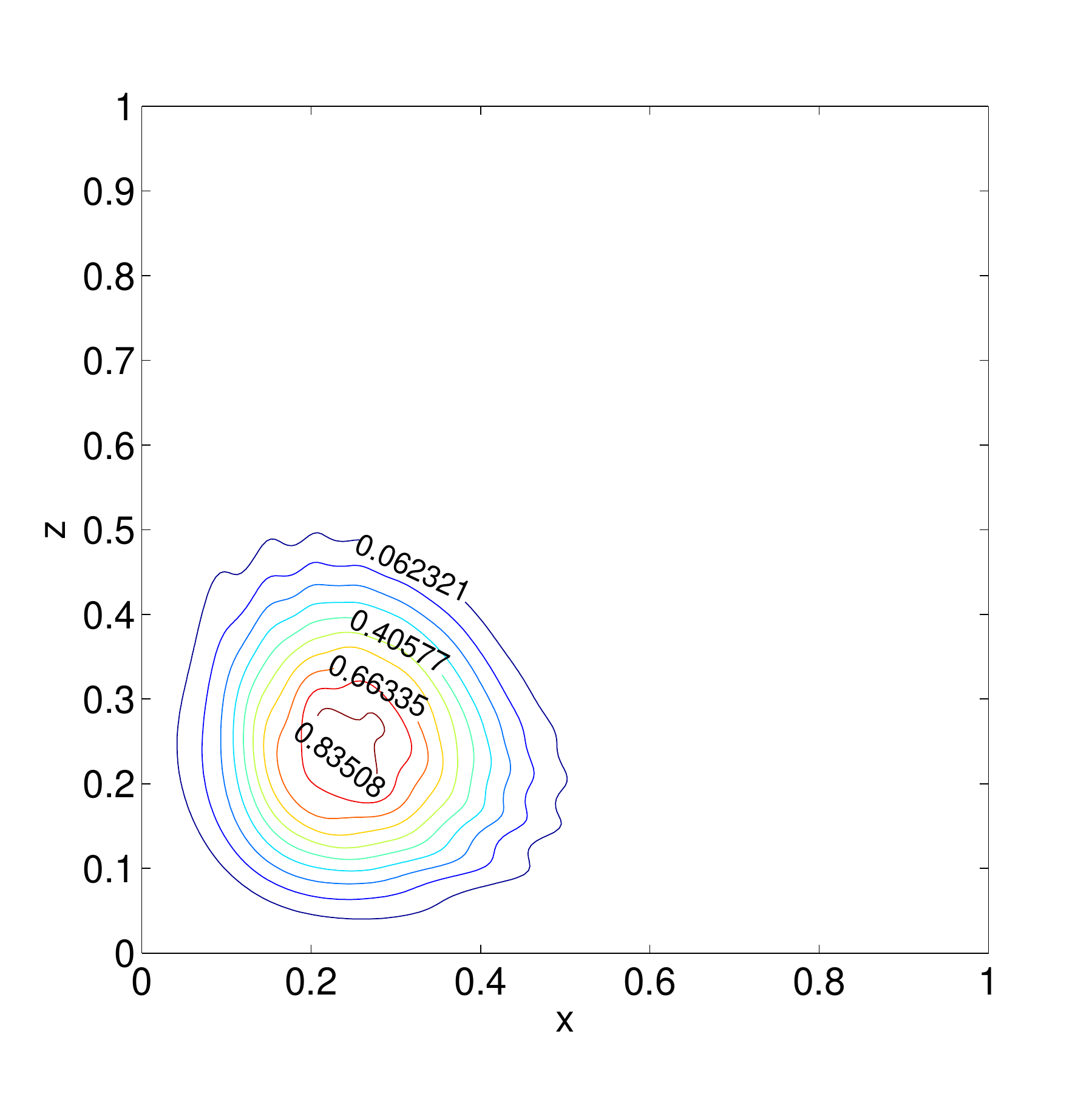}
\end{minipage}
\hspace{0.0cm}
\begin{minipage}[b]{0.5\linewidth}
\centering
\includegraphics[scale=0.3]{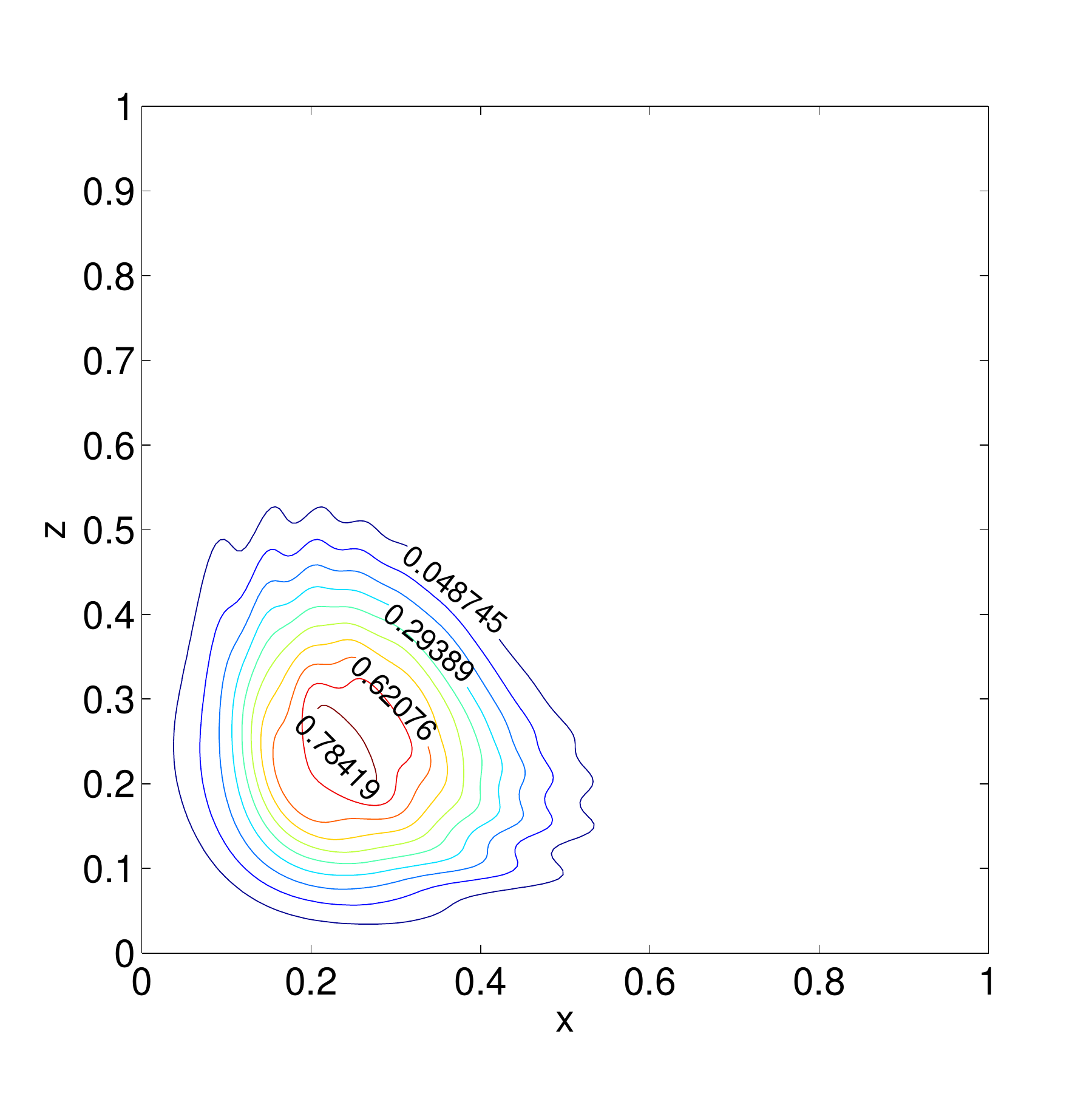}
\includegraphics[scale=0.3]{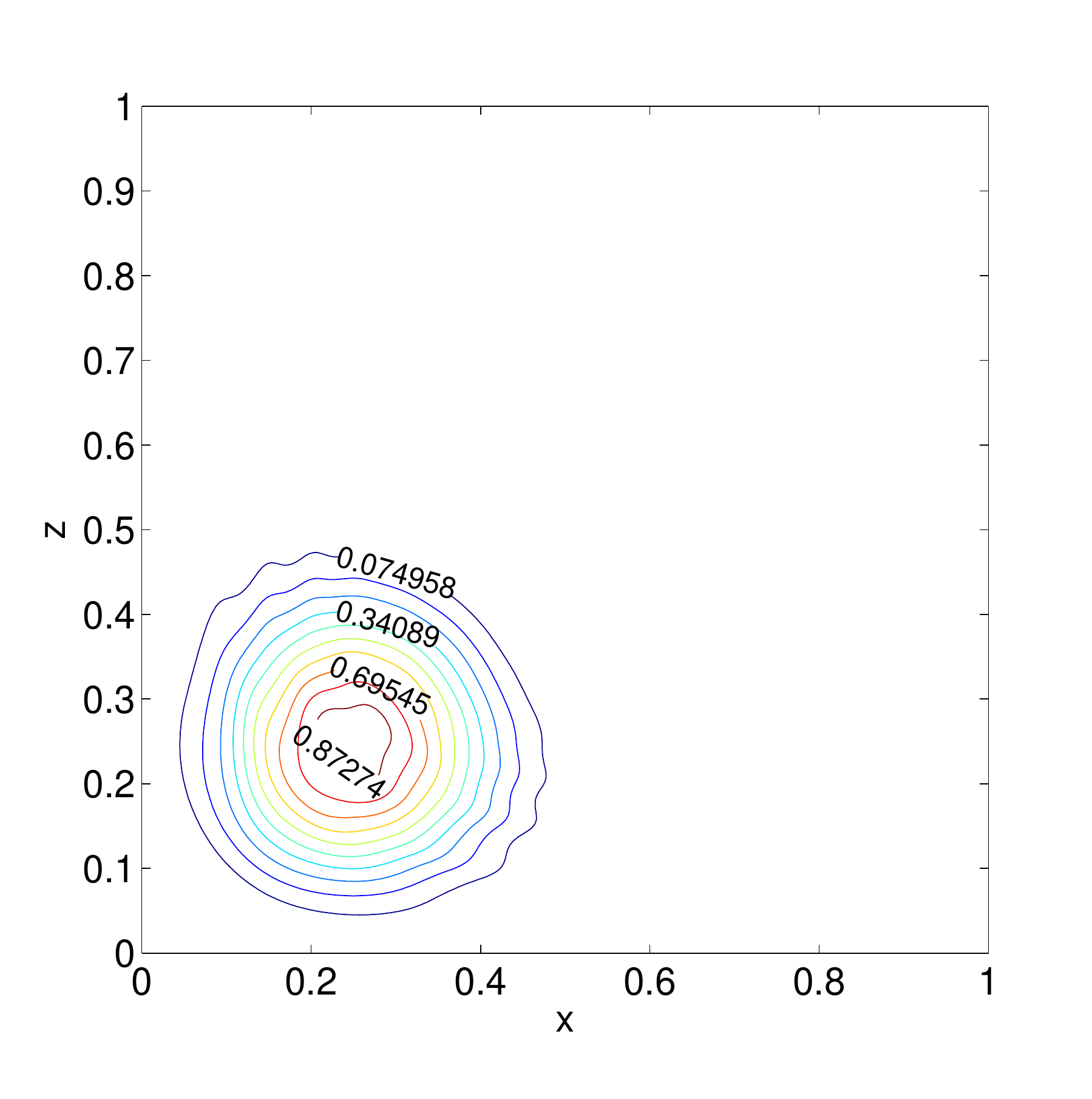}
\end{minipage}
\label{swirlinghigh}
\end{figure}

In this test we also include a study with low and high resolution simulations ($\dx=0.01$ and $\dx=0.005$ respectively), where we try different values for the weighting parameter in the GFORCE flux, $\omega$. Although the GFORCE flux used as low-order term in the FLIC approach is a first-order flux, (excepting for the case $\omega=1$ which corresponds to a Lax-Wendroff flux), increasing $\omega$ yields to a flux that has a behavior similar to a second-order flux; this is observed in figures \ref{swirlinglow} and \ref{swirlinghigh}, where noticeable increase on the accuracy of the final state is detected when switching to larger values of $\omega$ . Note that in every case the monotonicity of the solution is preserved. There is an additional cost related with the increase of accuracy and preserving monotonicity at the same time, which is a decrease on the time stepping. For instance, according to \cite{toroforce}, if $\omega$ is between 0.5 and 1,
\begin{equation}
\frac{\dt\max_{\Omega}|a|}{dx},\frac{\dt\max_{\Omega}|b|}{dz}\leq\left|\frac{-1+\omega}{2\omega}\right|,
\end{equation}
which implies that, for instance, when $\omega=0.75$, then $CFL\leq 0.17$.

\subsection{Doswell frontogenesis}
\begin{figure}
\caption{Doswell frontogenesis problem. Results at $t=4$[s] for smooth initial condition with $\delta=1$ and $200\times 200$ elements. In the right bottom, a cut at $x=0$.}
\centering
\includegraphics[width=\textwidth]{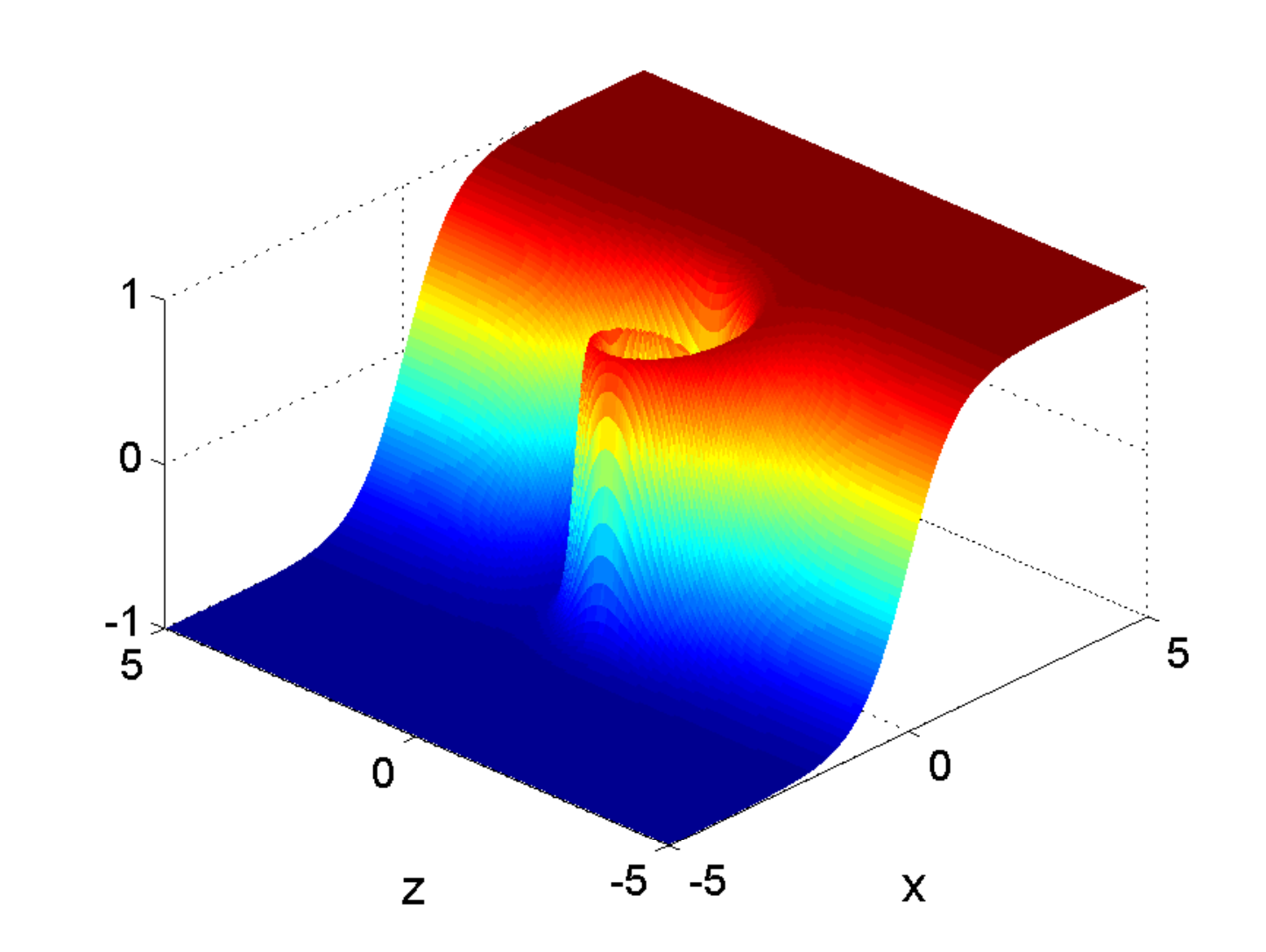}
\includegraphics[width=\textwidth, height=0.3\textheight]{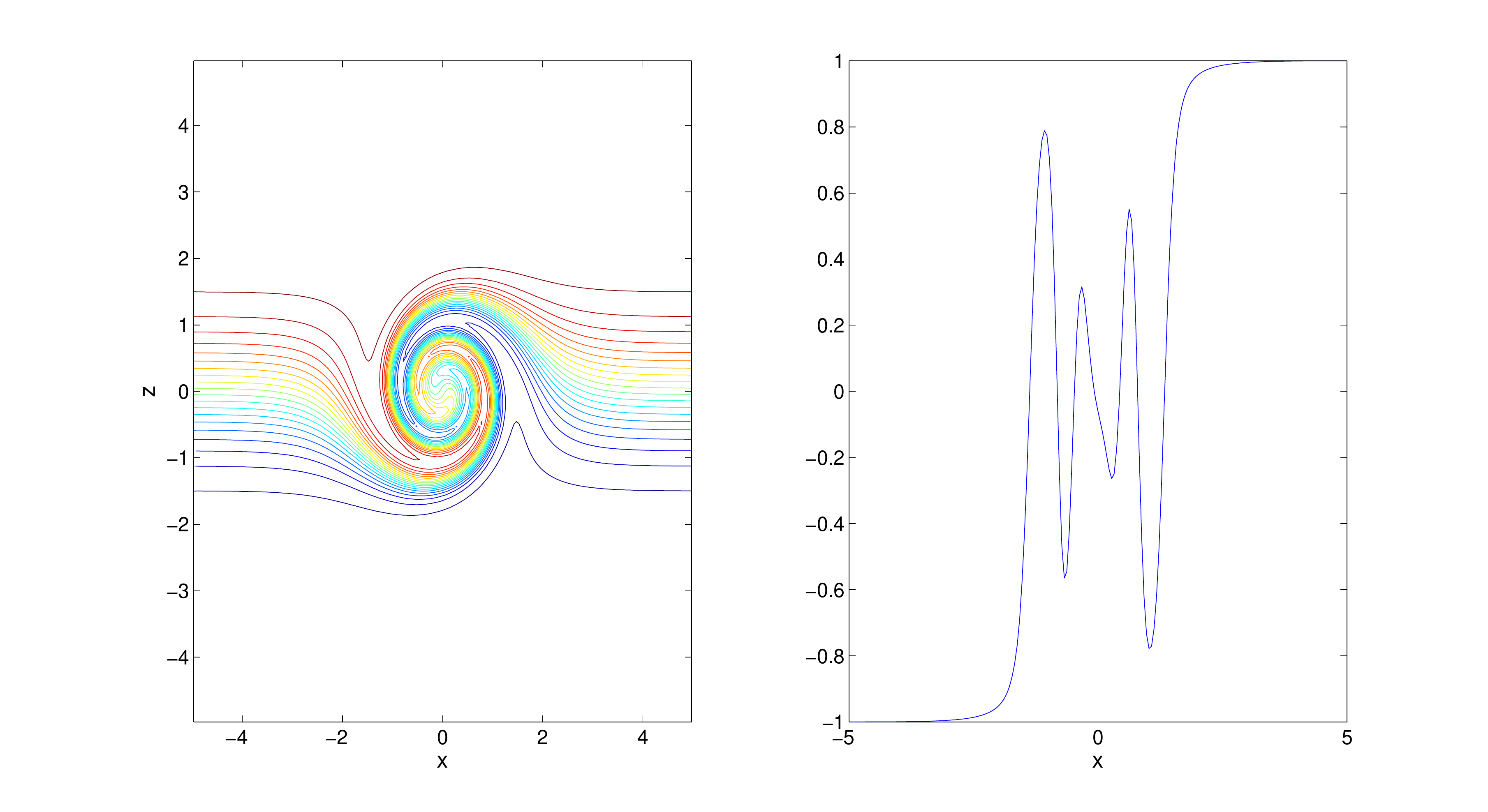}
\label{doswellsmooth}
\end{figure}

\begin{figure}
\caption{Doswell frontogenesis problem. Results at $t=4$[s] for sharp initial condition with $\delta=10^{-6}$ and $200\times 200$ elements. In the right bottom, a cut at $x=0$ with different limiters.}
\centering
\includegraphics[width=\textwidth]{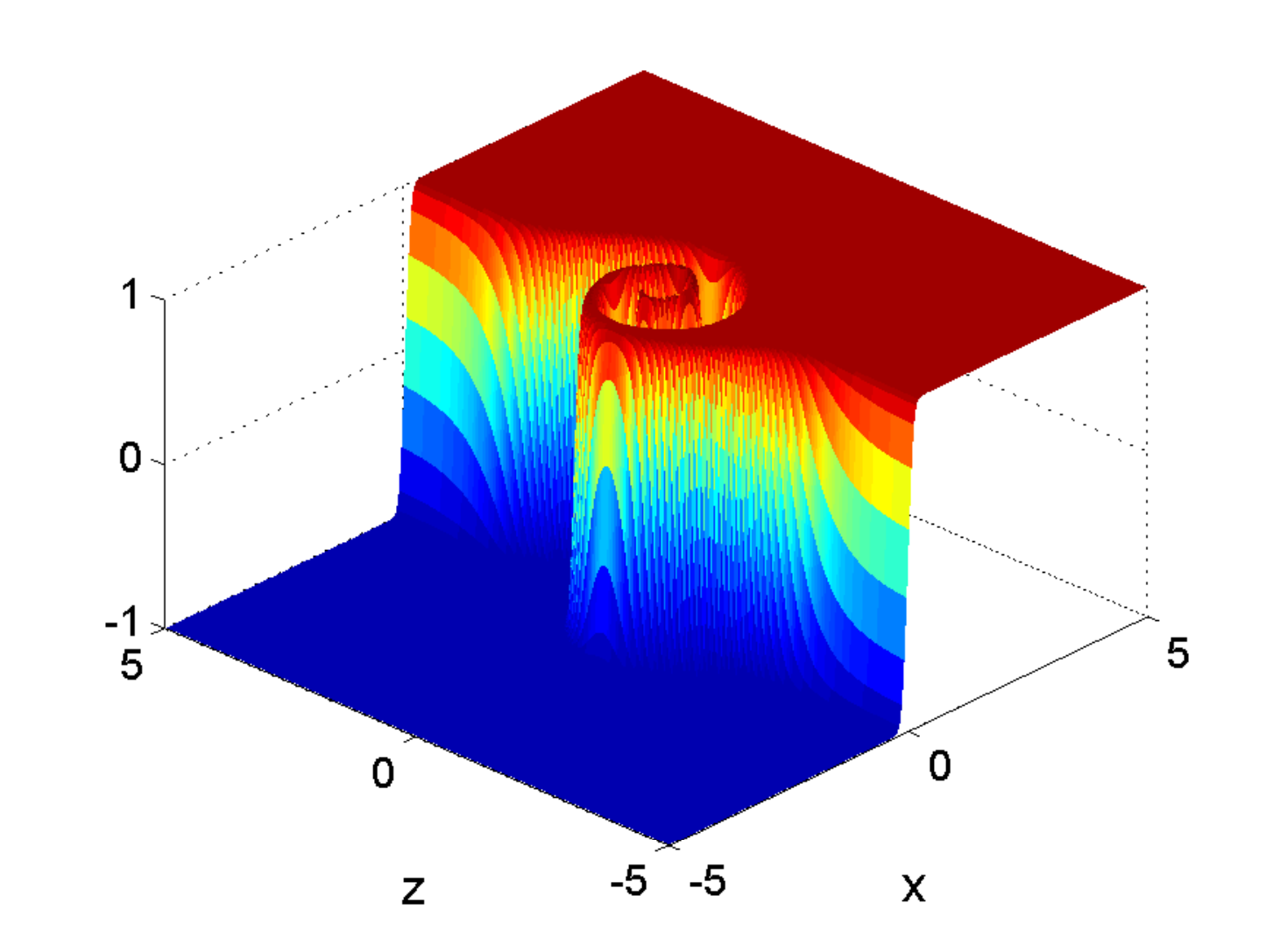}
\includegraphics[width=\textwidth, height=0.3\textheight]{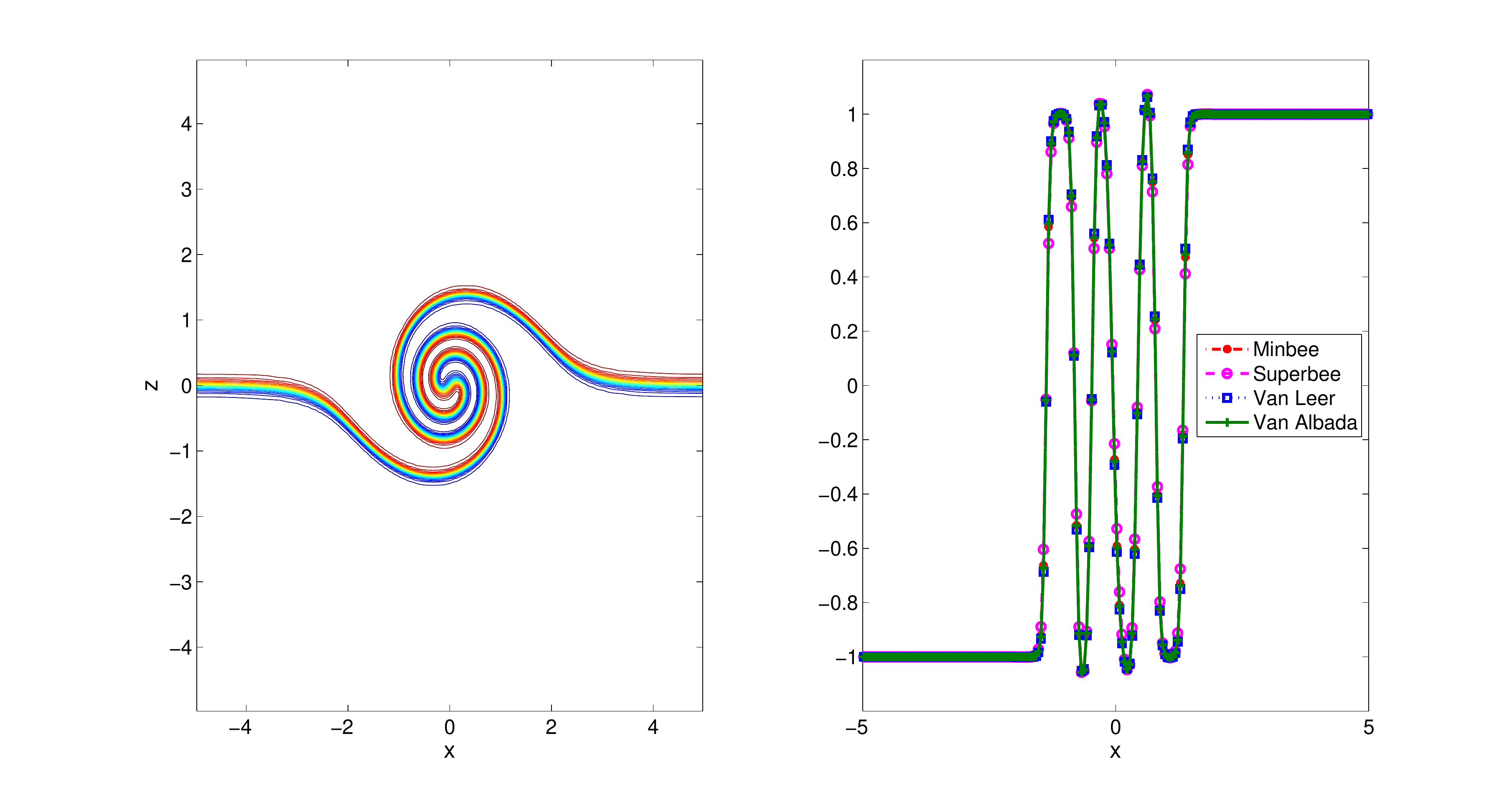}
\label{doswellsharp}
\end{figure}

We conclude the advection tests with a kinematic frontogenesis problem, originally presented in \cite{doswell}. It is a standard test in atmospheric modelling, and allows us to assess the performance of the scheme in the treatment of sharp fronts; numerical experiments with this test case in the context of this article can be found in \cite{titarefronto,gao}.
For this test, we set the domain $\Omega=[-5,\,5]^2$, and the flow is given by
\begin{equation}
a=-zf(r),\quad b=xf(r),\quad f(r)=\frac{1}{r}v(r),\quad v(r)=\bar v \,sech^2(r)\,tanh(r),
\end{equation}
\begin{equation}
r=\sqrt{x^2+z^2},\quad \bar v=2.59807.
\end{equation}
Initial condition for this test is given by
\begin{equation}
Q(x,z,0)=tanh\left(\frac{z}{\delta}\right),
\end{equation}
generating the following exact solution
\begin{equation}
Q(x,z,t)=tanh\left(\frac{z\cos(vt)-x\sin(vt)}{\delta}\right).
\end{equation}
The parameter $\delta$ is related with the thickness of the front zone. We first present numerical experiments with a value of $\delta=1$, in order to generate a smooth solution and verify the convergence rates; table \ref{tab:t3} shows that consistent convergence rates are obtained with the proposed scheme in both $L_{\infty}$ and $L_1$ norms. For this case, final state profiles can be observed in figure \ref{doswellsmooth}. A second is study is performed with sharp fronts, taking $\delta=10^{-6}$. Figure \ref{doswellsharp} illustrates the capapcity of the scheme in tracking sharp fronts; it can be observed that no spurious oscillations are generated, despite the sharpness of the solution. Moreover, the scheme proves to be very robust with respect to the choice of the limiter, which is interesting as a common situation in the choice of limiters is that it highly depends on the problem.

\begin{table}[H]
\centering
\caption{Convergence rates for the Doswell frontogenesis problem at  $t=4$ [s], with $\delta=1$.}\label{tab:t3}\vskip 3mm
\begin{tabular}{c c c c c}
\hline\\
N & $L_{\infty}$ error &  $L_{\infty}$ order & $L_1$ error &  $L_1$ order \\ [0.5ex]
\hline\\
50 &  3.4786e-001&  &1.2719e-002& \\
100 &   1.1140e-001 & 1.6 & 3.5136e-003 & 1.9  \\
200& 3.3302e-002 & 1.8 &  7.3045e-004 & 2.3 \\
400 &5.47780e-003 & 2.6  &  1.0541e-004 & 2.8  \\
\hline
\end{tabular}
\end{table}

\section{Convective tests: the Euler equations}
Throughout this section we consider a single model in order to study convective phenomena. Our starting point corresponds to the set of equations describing the evolution in time of a 2D dry air atmosphere \cite{kalnay,durran,torolibro}. Imposing conservation of mass, momentum and energy, and considering effects of gravity, together with neglecting friction and rotation effects, leads to a set of 2D inviscid primitive equations for the atmosphere written in conservative form
\begin{equation}\label{euler3d}
\pt Q + \px \F +\pz \H = \S,
\end{equation}
where
\begin{equation}
Q=\left[\begin{array}{c}\rho\\ \rho u\\\rho w\\ \rho\theta\end{array}\right],
\F=\left[\begin{array}{c}\rho u\\ \rho u^2 +\P\\ \rho uw\\ \rho u \theta\end{array}\right],
\H=\left[\begin{array}{c}\rho w\\ \rho wu\\\rho w^2+\P\\ \rho w \theta\end{array}\right],
\S=\left[\begin{array}{c}0\\0\\ -\rho g\\ 0\end{array}\right].
\end{equation}
In this system of equations $\rho$ is the density of the fluid, $u$ is the velocity in the x-direction, $w$ is the velocity in the z-direction, $\P$ is the pressure and $\theta$ is the potential temperature, which relates to the usual thermodynamic temperature $\mathcal{T}$ via
\begin{equation}
\theta\equiv \mathcal{T}\left(\frac{\P}{\P_0}\right)^{-R_d/c_p}.
\end{equation}
The system is closed by the equation of state for an ideal gas
\begin{equation}
\P=C_0(\rho\theta)^{\gamma},\qquad C_0=\frac{R_d^\gamma}{\P_0^{R_d/c_v}}.
\end{equation}
Model parameters are: the gravitational acceleration $g=9.81 [ms^{-2}]$, the atmospheric pressure at sea level $\P_0=10^5 [Pa]$, the gas constant for dry air $R_d=287 [JK^{-1}kg^{-1}]$, the specific heat of dry air at constant pressure and volume $C_p=1004 [JK^{-1}kg^{-1}]$, the specific heat of dry air at constant volume $C_v=717 [JK^{-1}kg^{-1}]$ and its ratio $\gamma=C_p/C_v=1.4$. Additionally, by defining the Exner pressure $\pi$
\begin{equation}
\pi\equiv\left(\frac{\P}{\P_0}\right)^{R/c_p}
\end{equation}
the expression for the total energy  of the system (internal+kinetic+potential) is given by
\begin{equation}\label{energia}
e=E_{Int}+E_{Kin}+E_{Pot}=c_v\theta\pi+\frac12(u^2+w^2)+gz,
\end{equation}
which will be used as flow parameter in the limiter computation. For all the simulations in this section, the parameter $\omega$ is set to 0.5, and once $\dx$ and $\dz$ (although always $\dx=\dz$) have been specified, the time stepping is selected according to
\begin{equation}
\dt=CFL\,\,\min\left(\frac{\dx}{\displaystyle{\max_{\Omega}}\,s_x}, \frac{\dz}{\displaystyle{\max_{\Omega}}\,s_z}\right),
\end{equation}
with $CFL=0.4$ and where $s_x$ and $s_z$ are maximum characteristic speeds in the $x$ and $z$ direction respectively,
\begin{equation}
s_x=\max (u+c_s,u-c_s),\quad s_z=\max (w+c_s,w-c_s),
\end{equation}
where $c_s=\sqrt{\partial_{\rho}\P}$ is the speed of sound in the fluid. In this section we study four test cases for the above presented set of equations. They all consist of hydrostatically balanced initial conditions plus a potential temperature perturbation. The first three cases are initialized with reference states for a neutral atmosphere, i.e.,

\begin{equation}
\theta=\bar\theta+\theta',\qquad\bar\theta=300[K],
\end{equation}
and the density $\rho$ is initialized via the hydrostatic balance equation
\begin{equation}
c_p\theta\frac{d\pi}{dz}=-g,
\end{equation}
together with the ideal gas law
\begin{equation}
\rho=\frac{P_0}{R_d\theta}\pi^{\frac{c_v}{R_d}},
\end{equation}
both evaluated at the reference state $\theta=\bar\theta$. The potential temperature perturbation, together with the initial velocity field and boundary conditions are specified for each problem.

\subsection{Convective bubble in a neutral atmosphere}

This first test case, that has been previously addressed in \cite{robert,smolarbubble,giraldo} (among others), studies the behavior of a hot temperature bubble placed in a hydrostatically balanced neutral atmosphere at rest. As the perturbation is warmer than the background state, a buoyancy force will push it upwards and, as it starts rising, because of the same buoyancy effect, it will start experiencing a deformation that will eventually develop into a mushroom-type of cloud. We use the scale and settings used in \cite{knothjcp}. The domain is $\Omega=[-10000,\,10000]\times[0,\,10000]$, with $\dx=\dz=125[m]$
and the potential temperature perturbation is given by
\begin{equation}
\theta'=\begin{cases}2\cos\left(\frac{\pi L}{2}\right) & L\leq 1,\\ 0 & i.o.c.\end{cases},\quad L=\frac{1}{2000}\sqrt{x^2+(z-2000)^2}.
\end{equation}
Simulation time has been set to 1000 [s], allowing the bubble to rise without hitting the top boundary; reflecting solid wall boundary conditions have been considered around the whole domain.

\begin{figure}
\caption{Convective bubble in a neutral atmosphere test case; potential temperature colormaps. From left to right, top to bottom: initial state and Results at $t=300$, 600 and 1000 [s], with $\Delta x=\Delta z=125$[m], $160\times80$ elements.}
\begin{minipage}[b]{0.5\linewidth}
\centering
\includegraphics[scale=0.5]{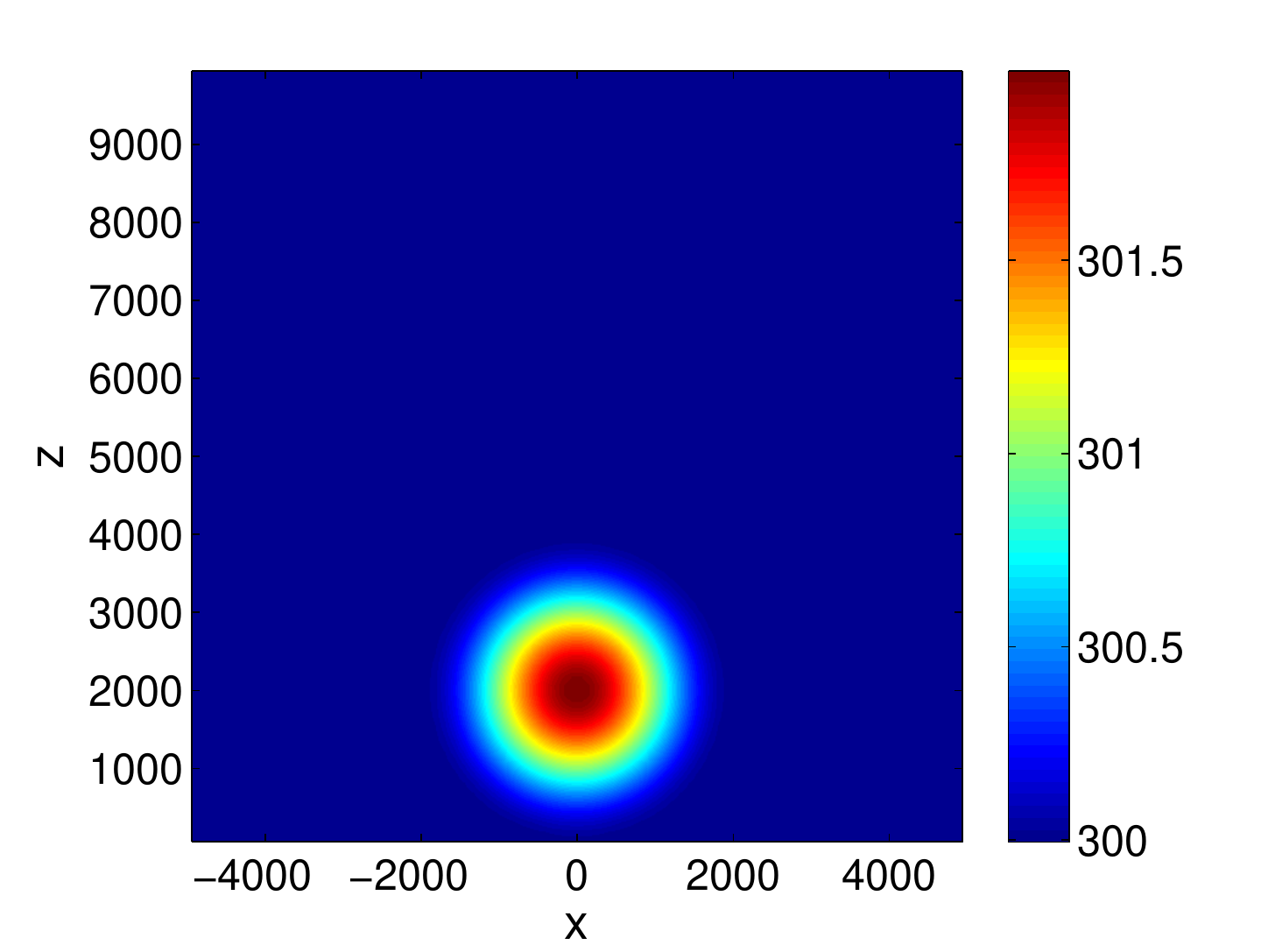}
\includegraphics[scale=0.5]{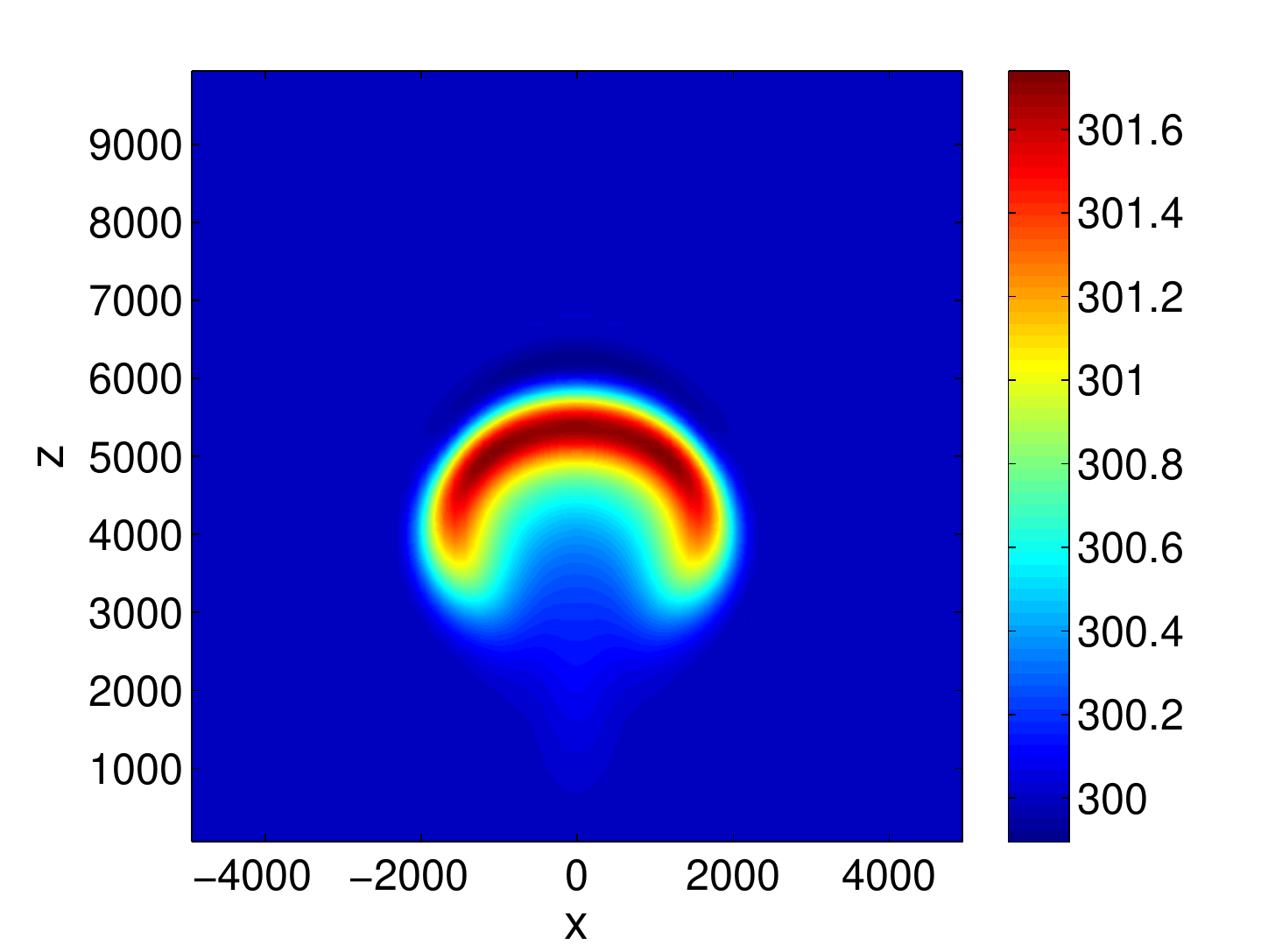}
\end{minipage}
\hspace{0.0cm}
\begin{minipage}[b]{0.5\linewidth}
\centering
\includegraphics[scale=0.5]{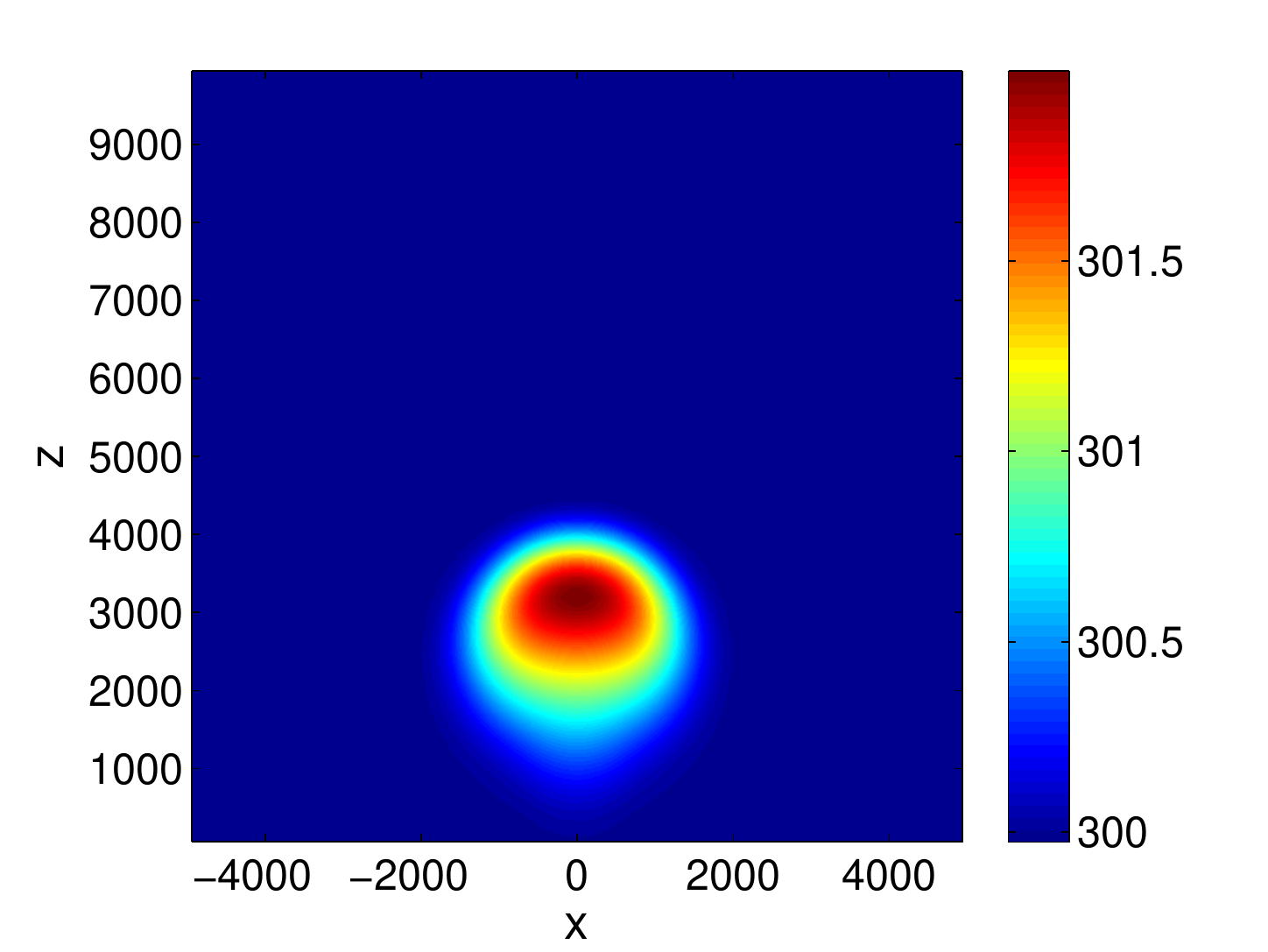}
\includegraphics[scale=0.5]{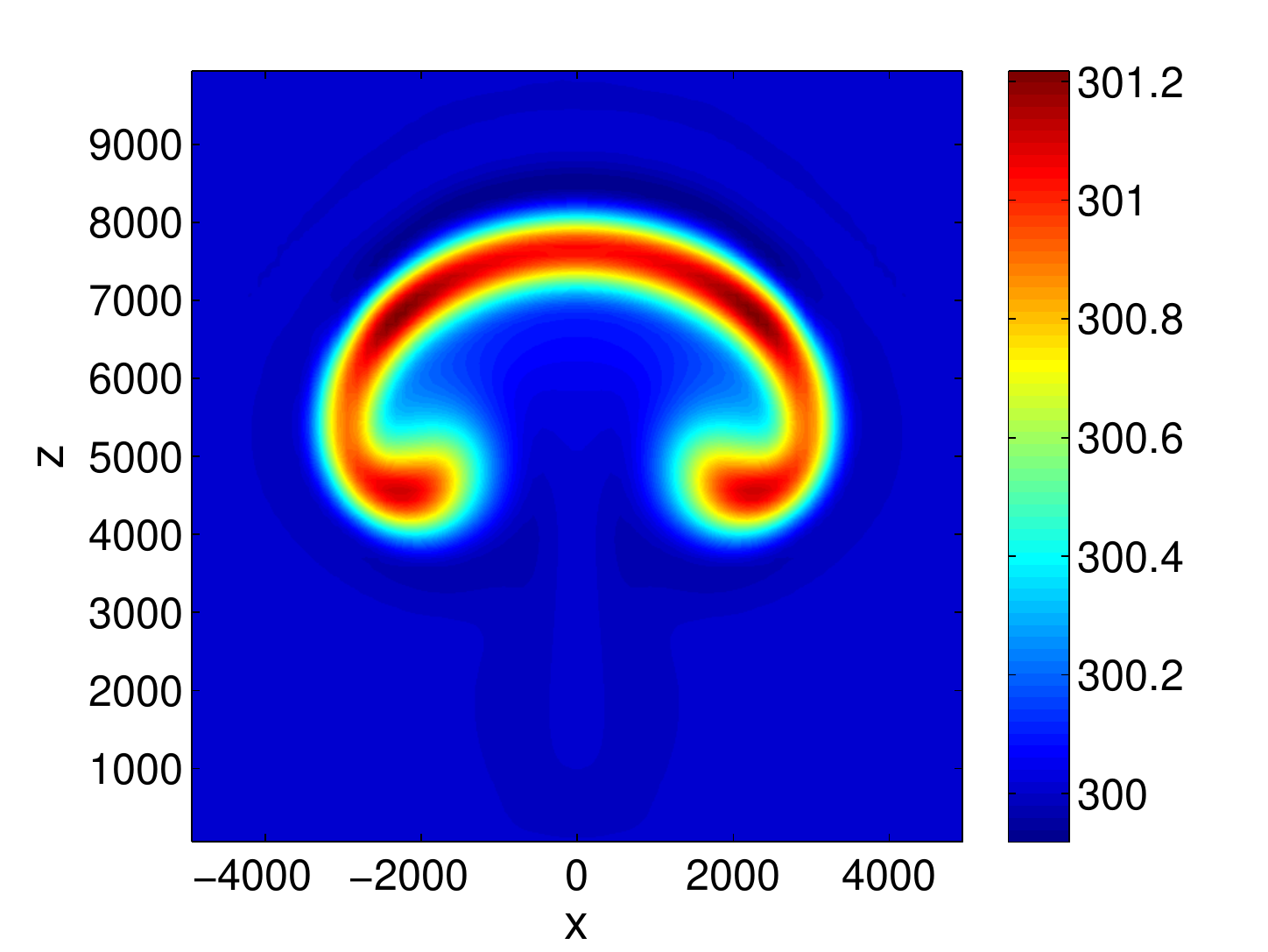}
\end{minipage}
\label{bubblet4}
\end{figure}

\begin{figure}
\caption{Convective bubble in a neutral atmosphere test case; velocity field colormaps at $t=600$ and $t=1000$ [s], with $\Delta x=\Delta z=125$[m], $160\times80$ elements. Left: horizontal velocity $u$. Right: vertical velocity $w$.}
\begin{minipage}[b]{0.5\linewidth}
\centering
\includegraphics[scale=0.5]{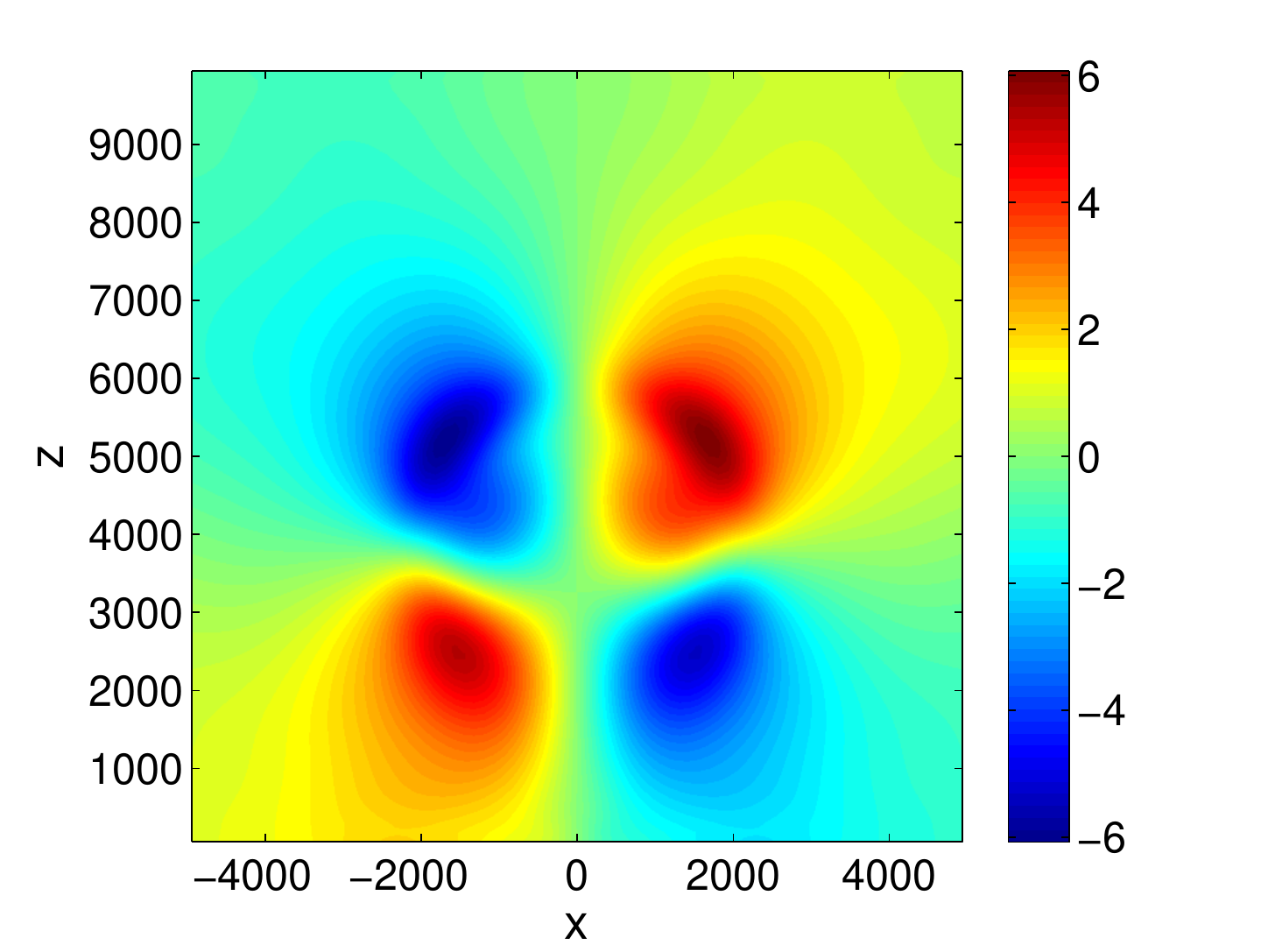}
\includegraphics[scale=0.5]{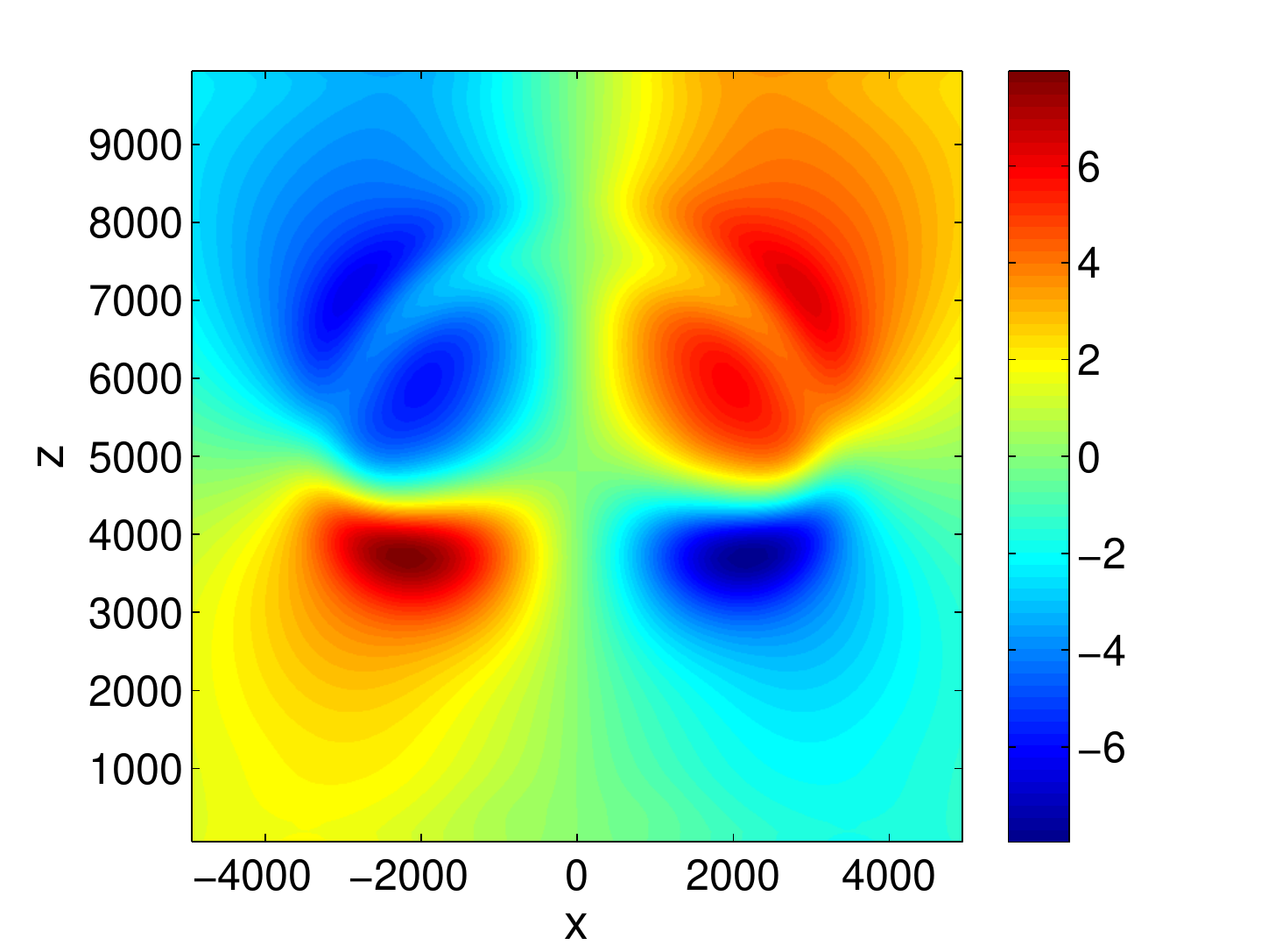}
\end{minipage}
\hspace{0.0cm}
\begin{minipage}[b]{0.5\linewidth}
\centering
\includegraphics[scale=0.5]{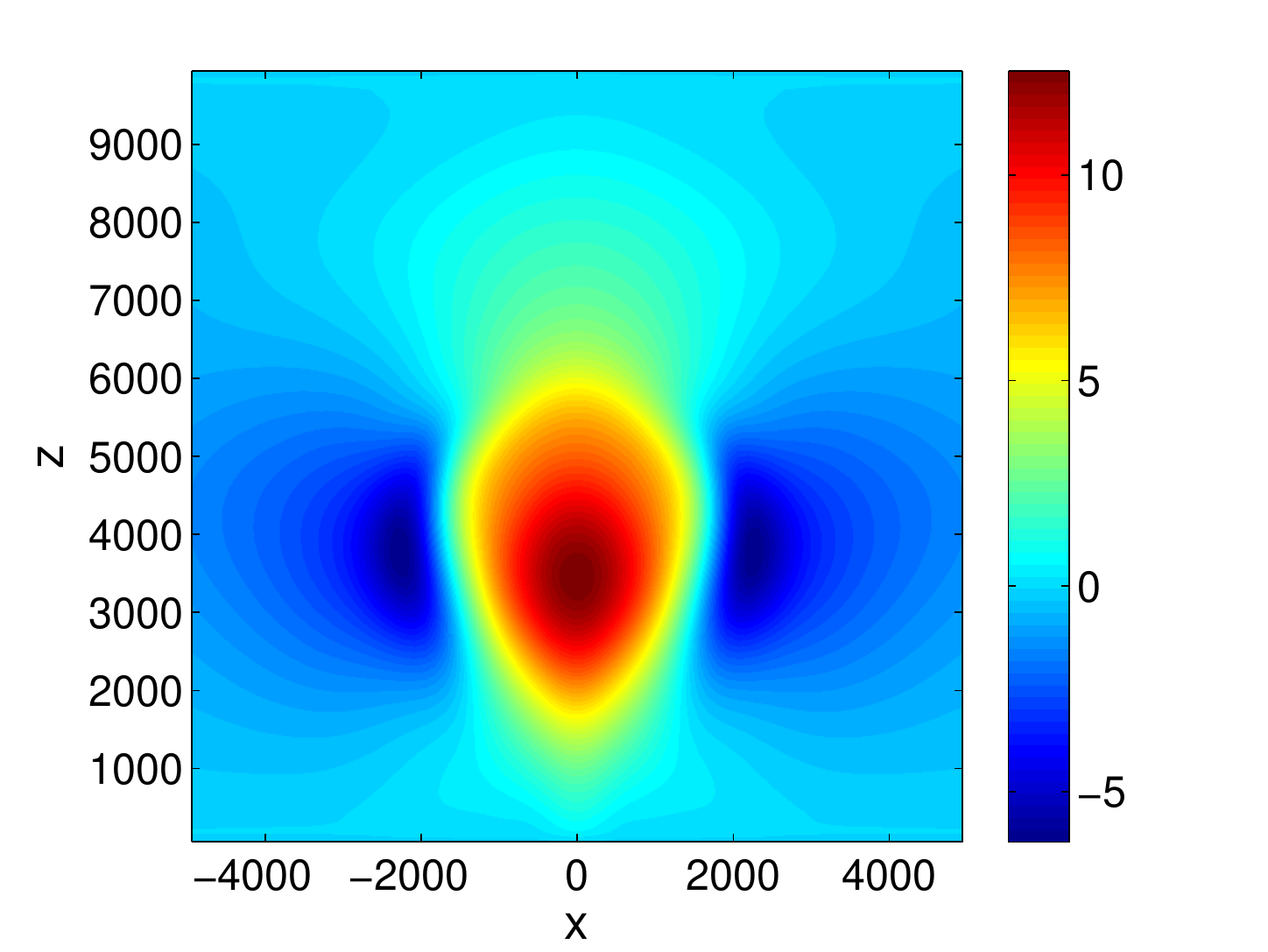}
\includegraphics[scale=0.5]{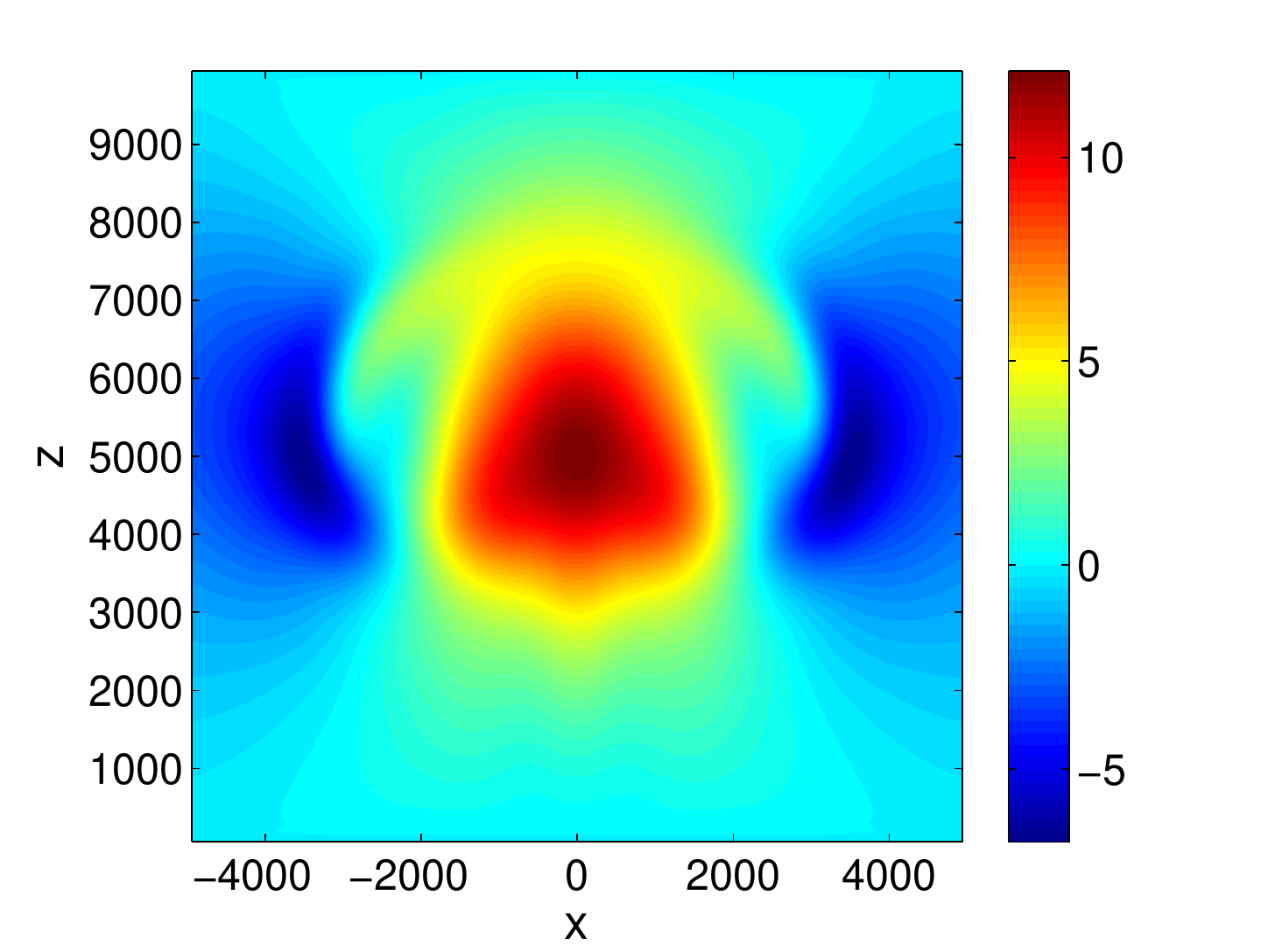}
\end{minipage}
\label{bubblev4}
\end{figure}

\begin{figure}
\caption{Normalized energy (with respect to the initial total value) of the system for the convective bubble in a neutral atmosphere test case.}
\centering
\includegraphics[height=0.3\textheight, width=\textwidth]{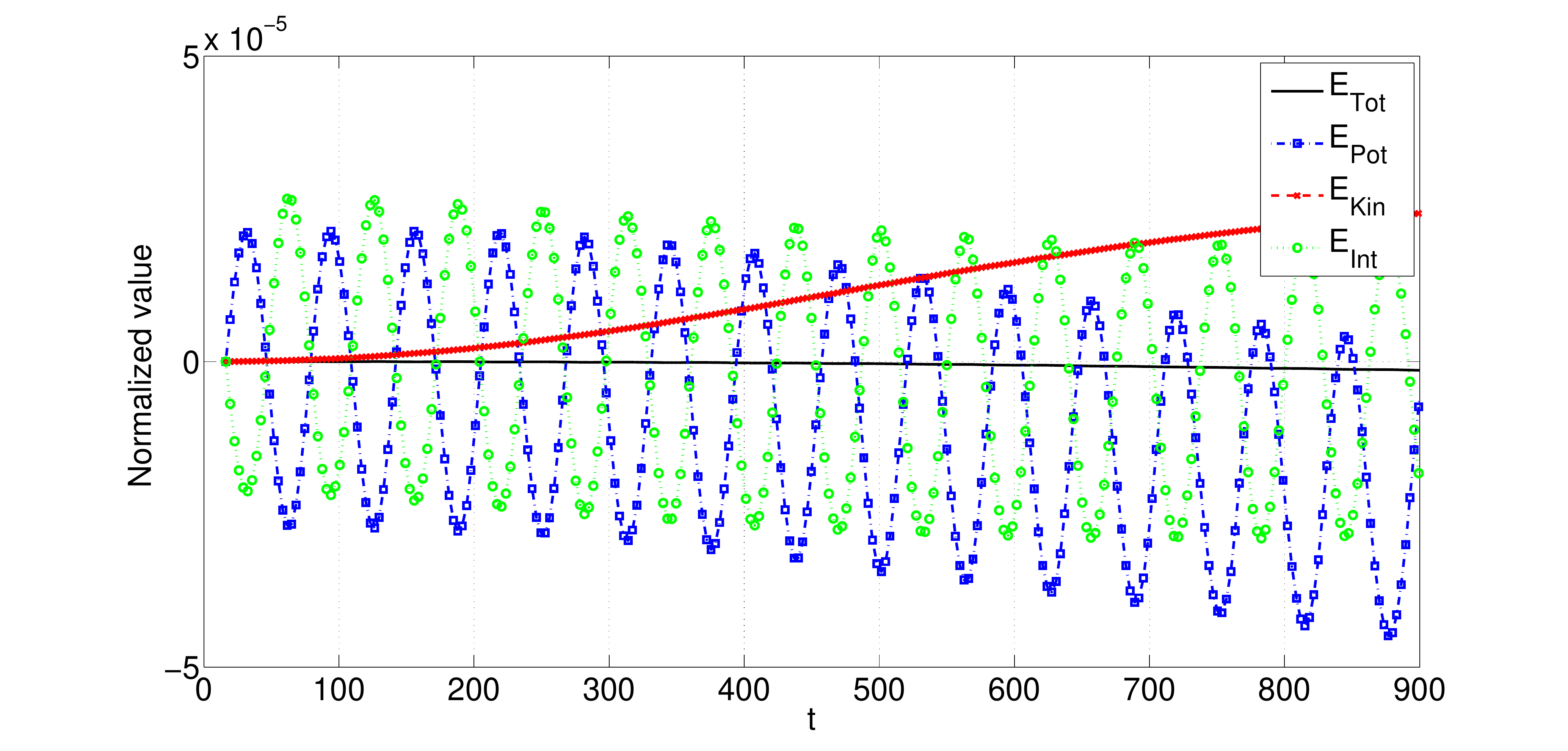}
\label{bubbleet4}
\end{figure}

It can be seen in figure \ref{bubblet4} that the proposed scheme manages to reproduce the correct physical solution of the test; the bubbles rises from its initial state experiencing deformation that finally generates a mushroom cloud. Note that symmetry is preserved during the simulation, no spurious oscillations are detected and the final front location is essentially the same as shown in \cite{knothjcp}, also matching in the vertical velocity, as shown in figure \ref{bubblev4}. It has to be remarked that most of the numerical models used to approximate such experiment need to include a damping effect in order to obtain a grid converged solution, which is not our case. We also assess the performance of the scheme in this test by the energy plot in figure \ref{bubbleet4}; it illustrates the quantities present in eq. (\ref{energia}) (we consider the the quantity $\rho e$ rather than $e$), and it can be seen the rise of the kinetic energy of the system initially at rest together with a decrease of the internal and potential energy contributions to the system, all this with a constant total energy, as it is expected for a closed system without physical dissipation mechanisms.

\subsection{Interaction between hot and cold bubbles in a neutral atmosphere}
\begin{figure}
\caption{Colormap of the initial potential temperature for the hot and cold bubbles test case.}
\centering
\includegraphics[scale=0.5]{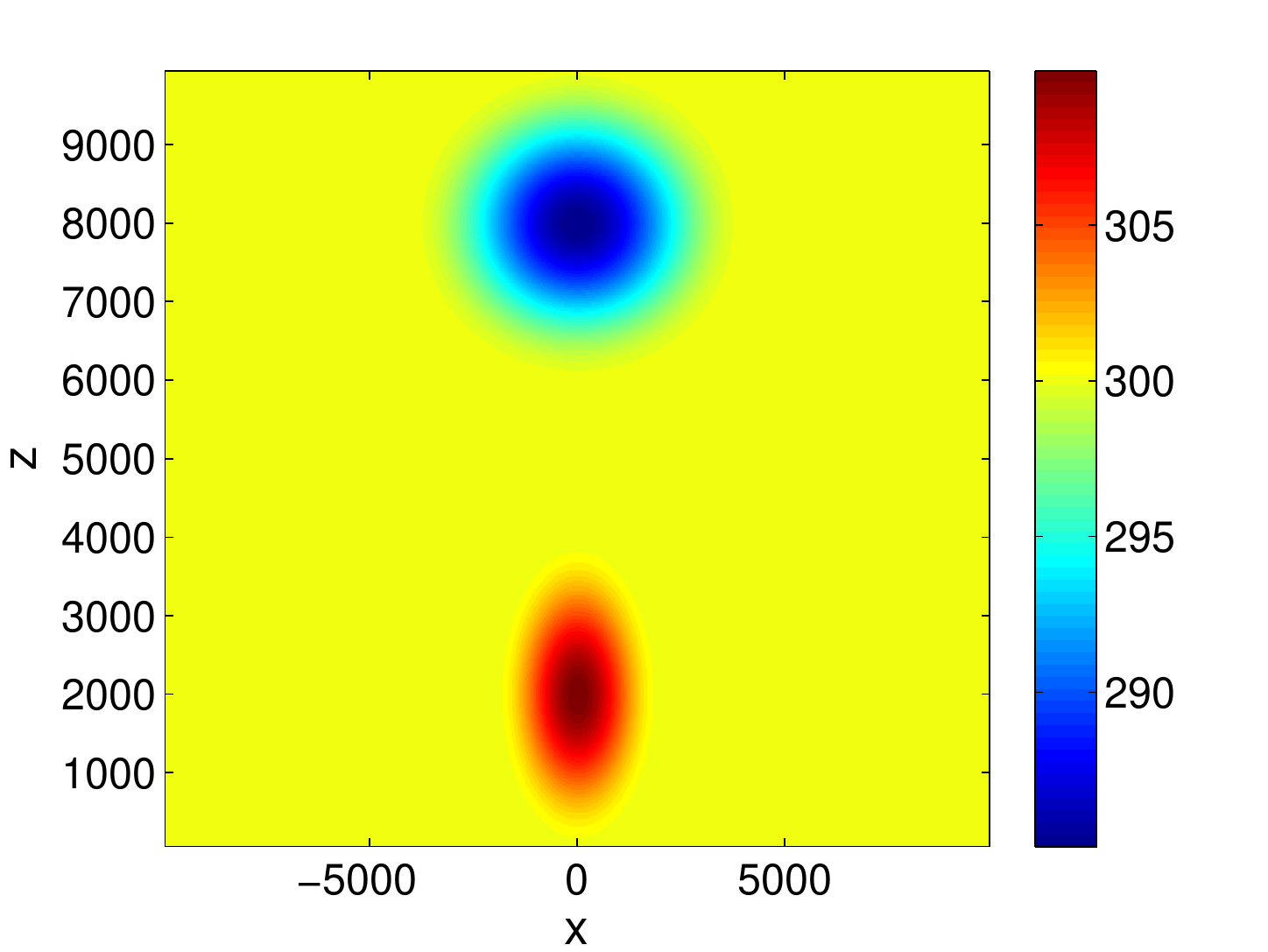}
\label{initt2}
\end{figure}

\begin{figure}
\caption{Hot and cold bubbles test case.Results at $t=180$ and 250[s], with $\Delta x=\Delta z=125$[m], $160\times 80$ elements. Left: colormap of the potential temperature. Right: vector plot of the velocity field.}
\begin{minipage}[b]{0.5\linewidth}
\centering
\includegraphics[scale=0.5]{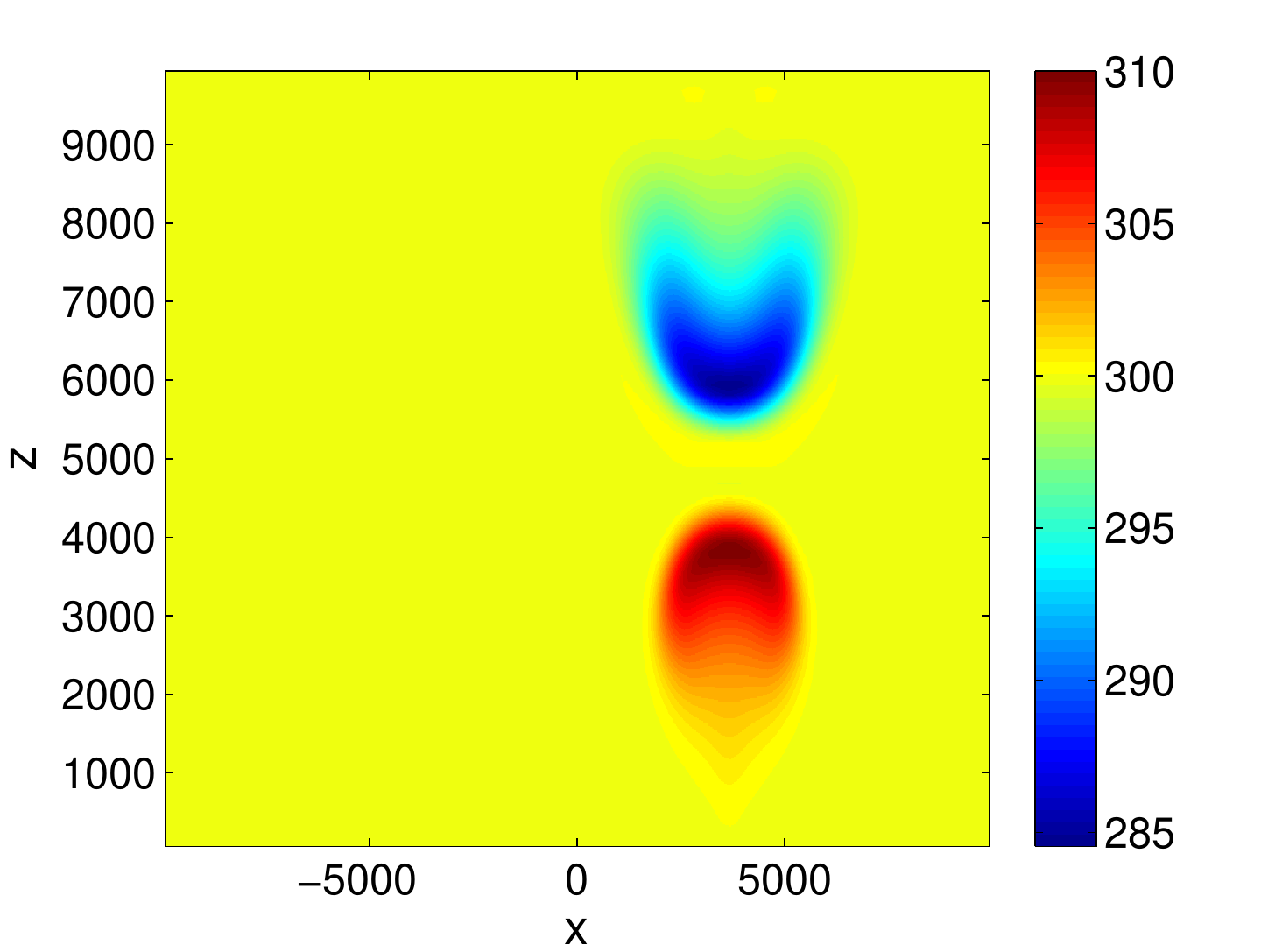}
\end{minipage}
\hspace{0.0cm}
\begin{minipage}[b]{0.5\linewidth}
\centering
\includegraphics[scale=0.5]{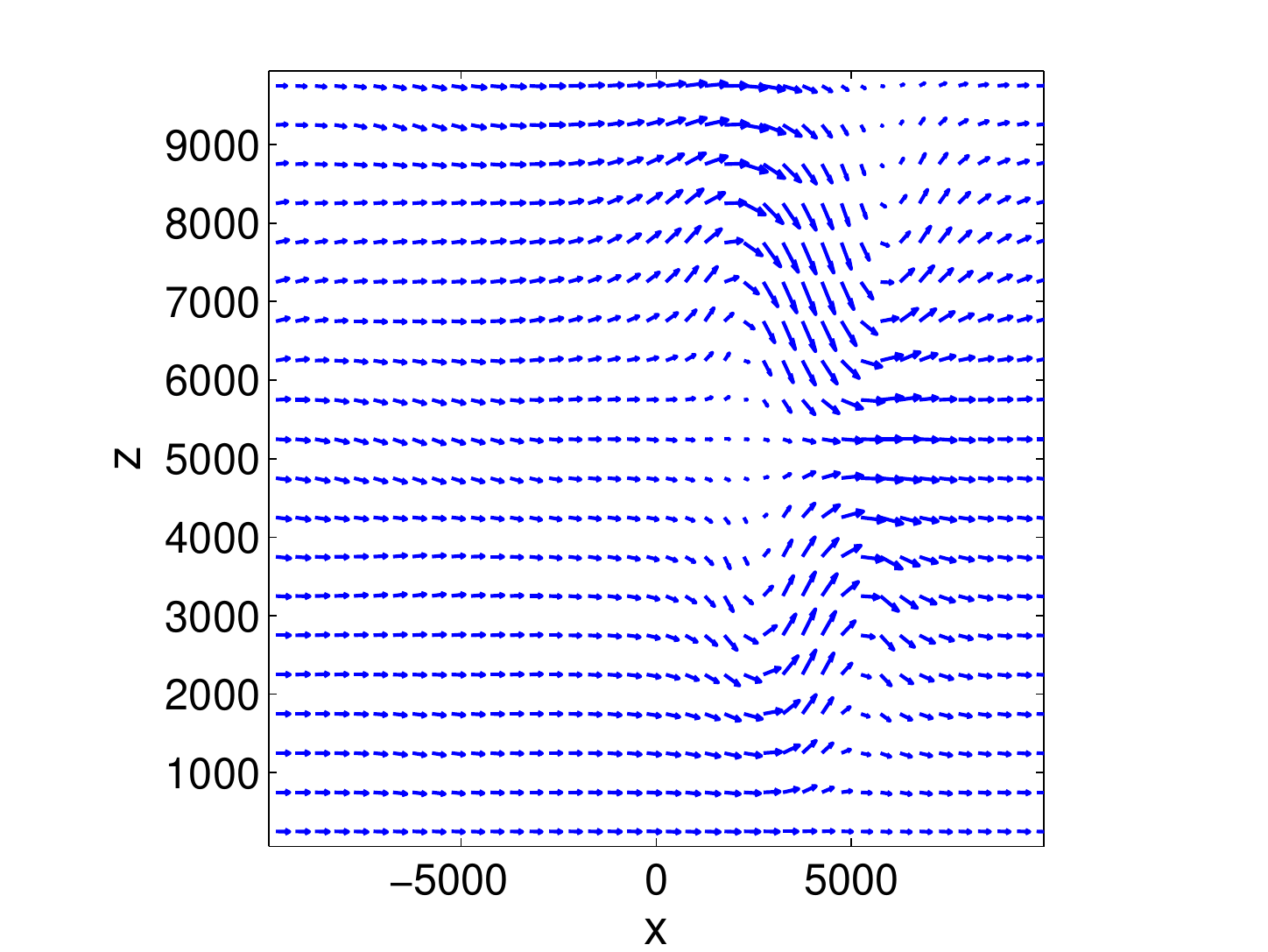}
\end{minipage}

\begin{minipage}[b]{0.5\linewidth}
\centering
\includegraphics[scale=0.5]{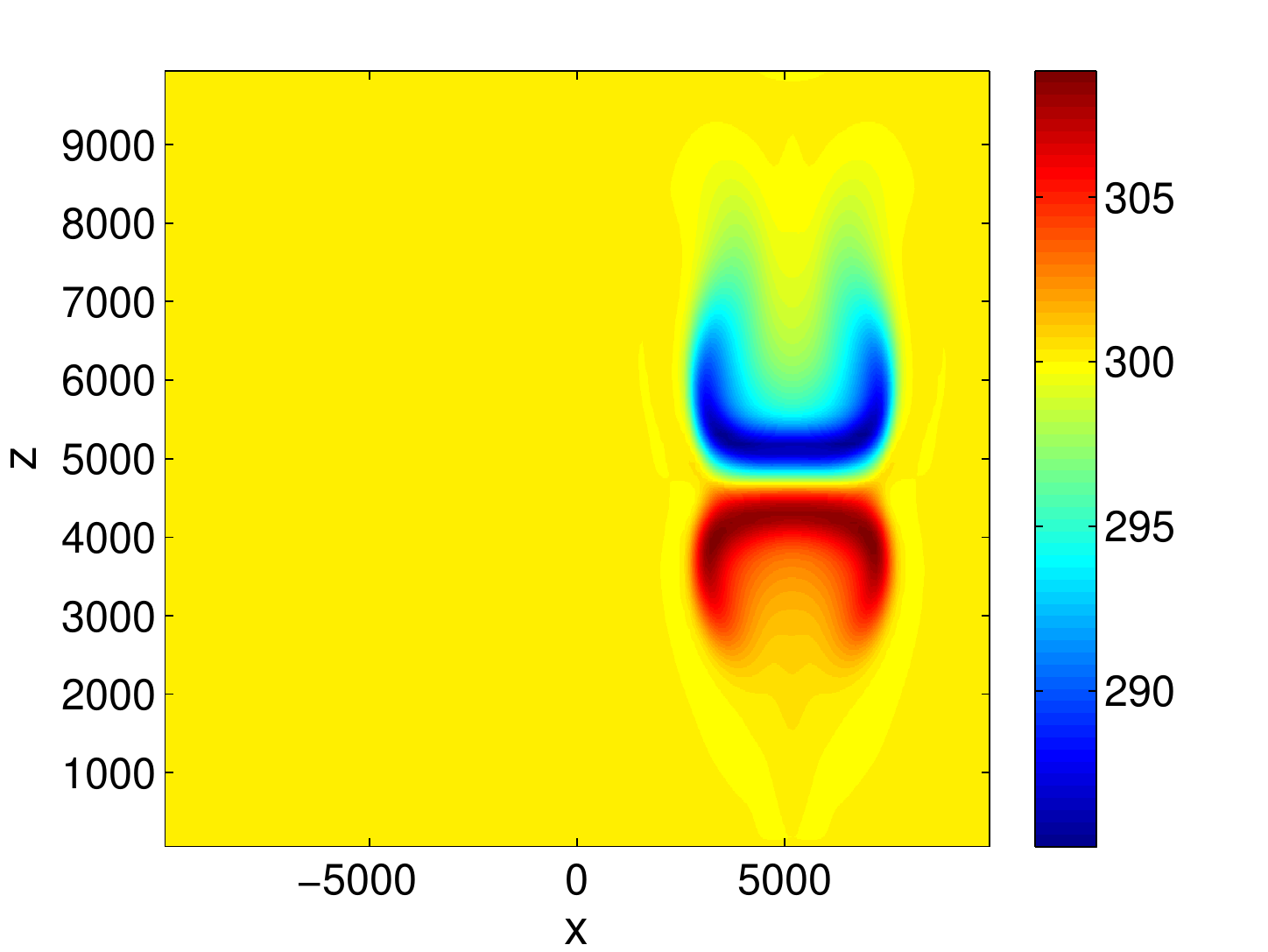}
\end{minipage}
\hspace{0.0cm}
\begin{minipage}[b]{0.5\linewidth}
\centering
\includegraphics[scale=0.5]{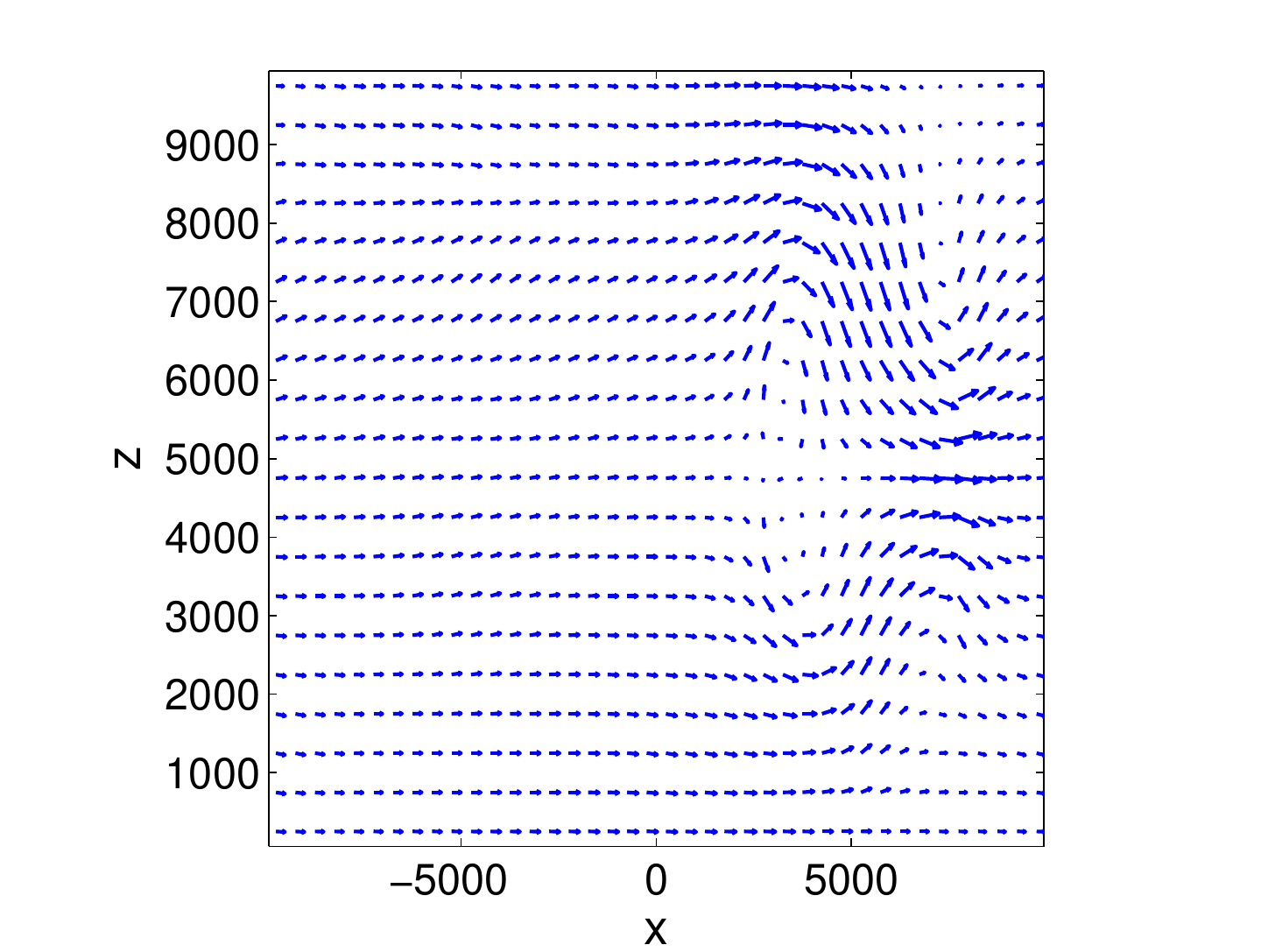}
\end{minipage}
\label{t2}
\end{figure}

\begin{figure}
\caption{Hot and cold bubbles test case.Results at $t=500$ and 1000[s], with $\Delta x=\Delta z=125$[m], $160\times80$ elements. Left: colormap of the potential temperature. Right: vector plot of the velocity field.}
\begin{minipage}[b]{0.5\linewidth}
\centering
\includegraphics[scale=0.5]{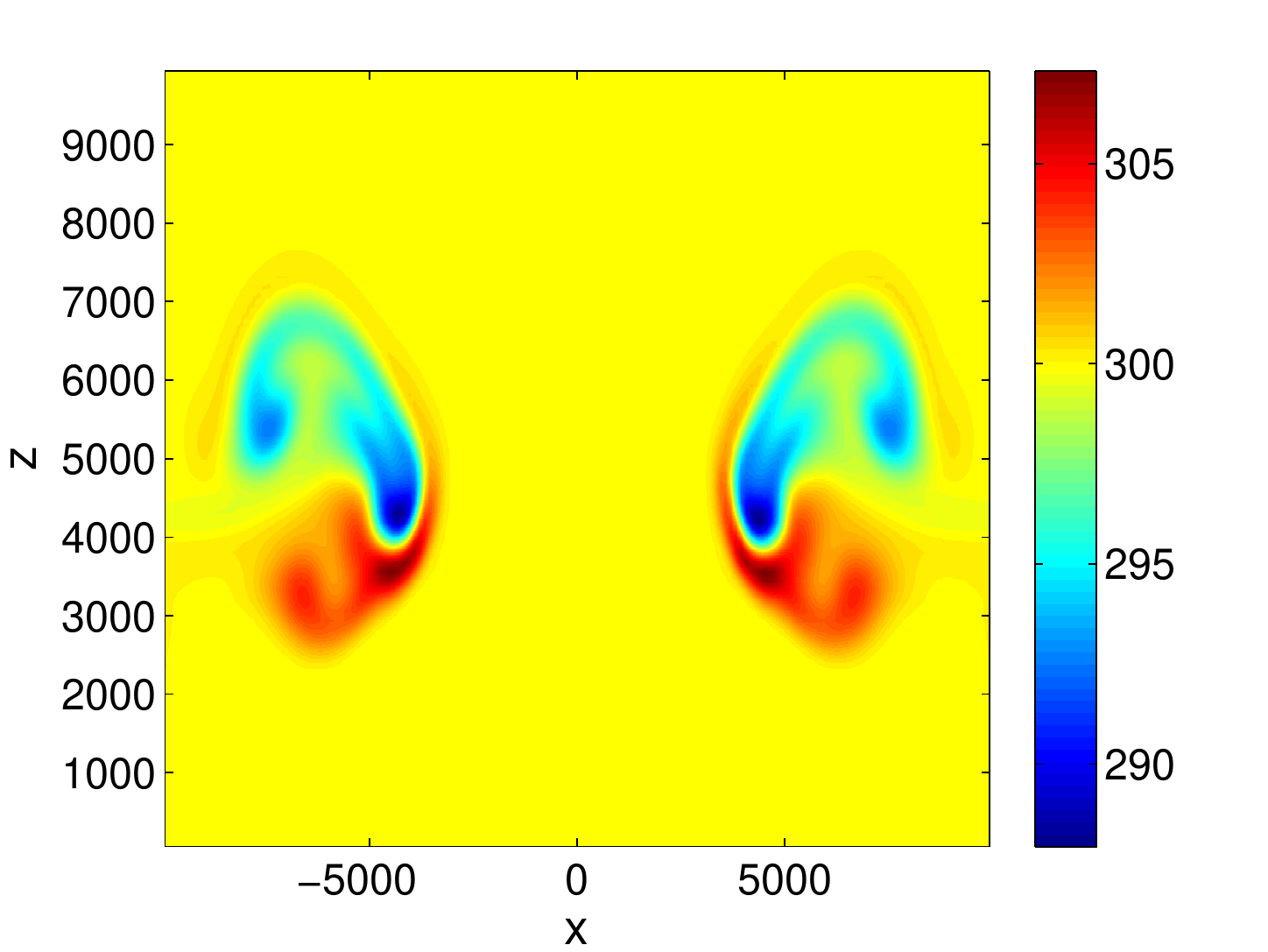}
\end{minipage}
\hspace{0.0cm}
\begin{minipage}[b]{0.5\linewidth}
\centering
\includegraphics[scale=0.5]{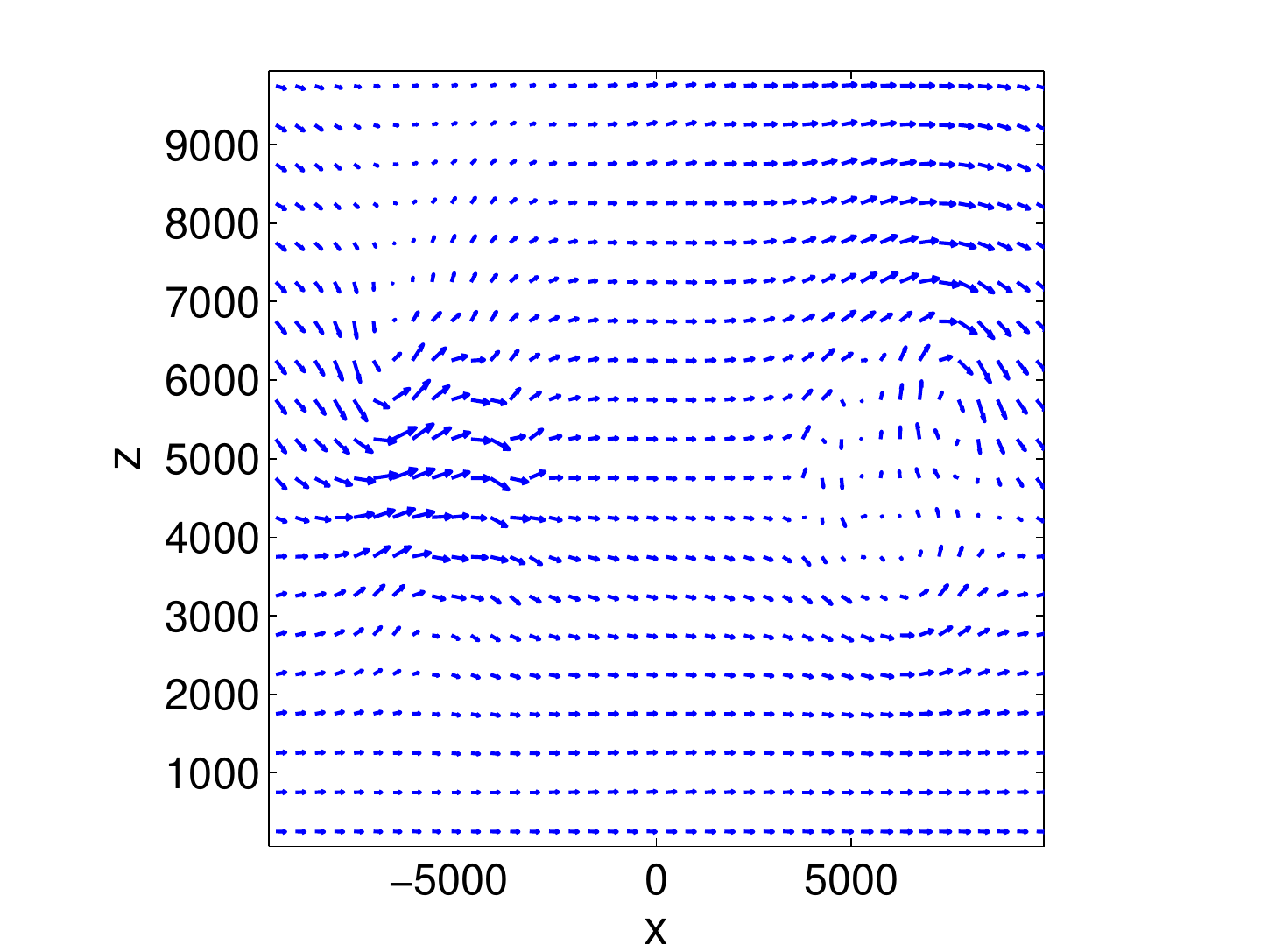}
\end{minipage}

\begin{minipage}[b]{0.5\linewidth}
\centering
\includegraphics[scale=0.5]{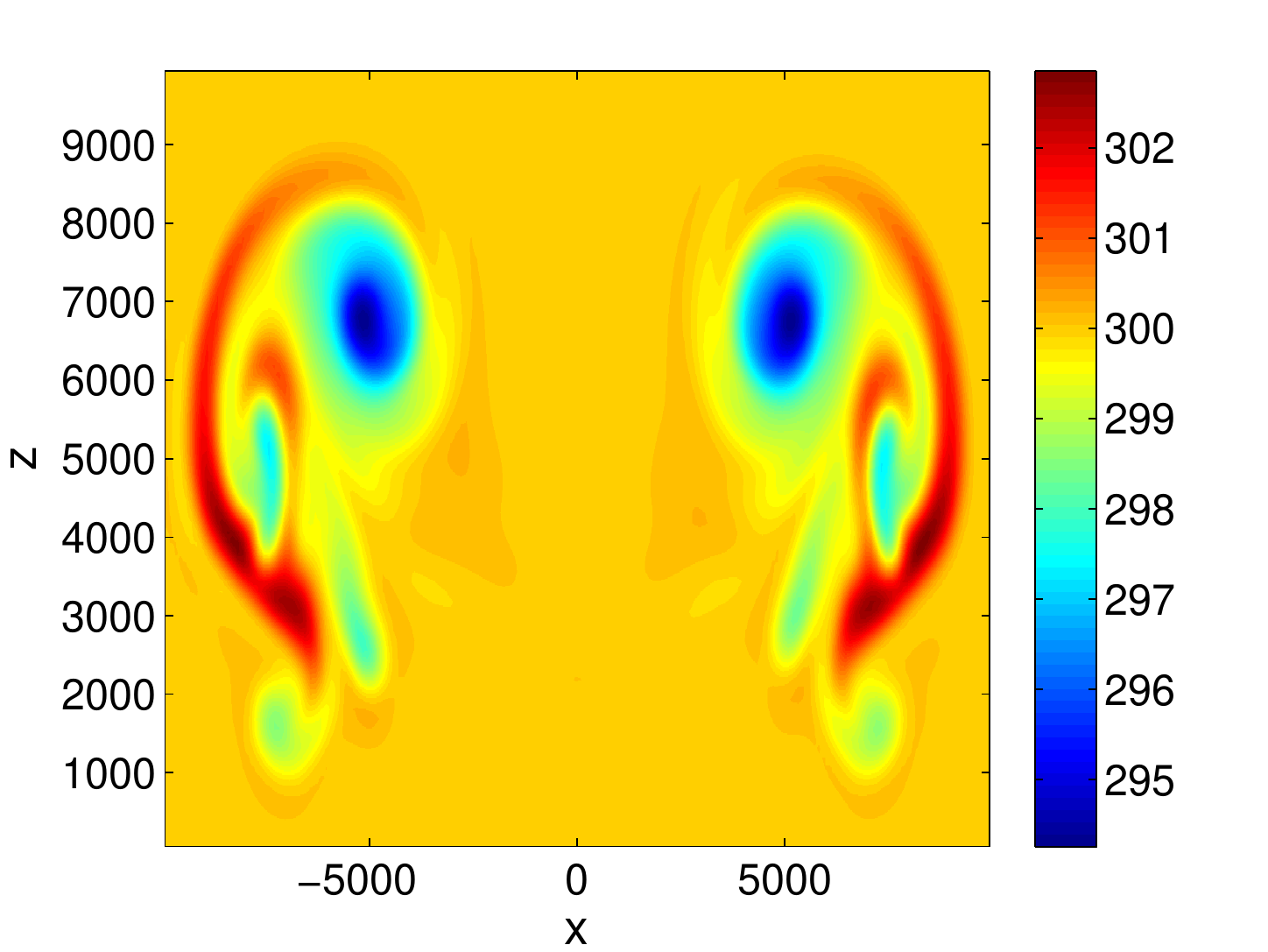}
\end{minipage}
\hspace{0.0cm}
\begin{minipage}[b]{0.5\linewidth}
\centering
\includegraphics[scale=0.5]{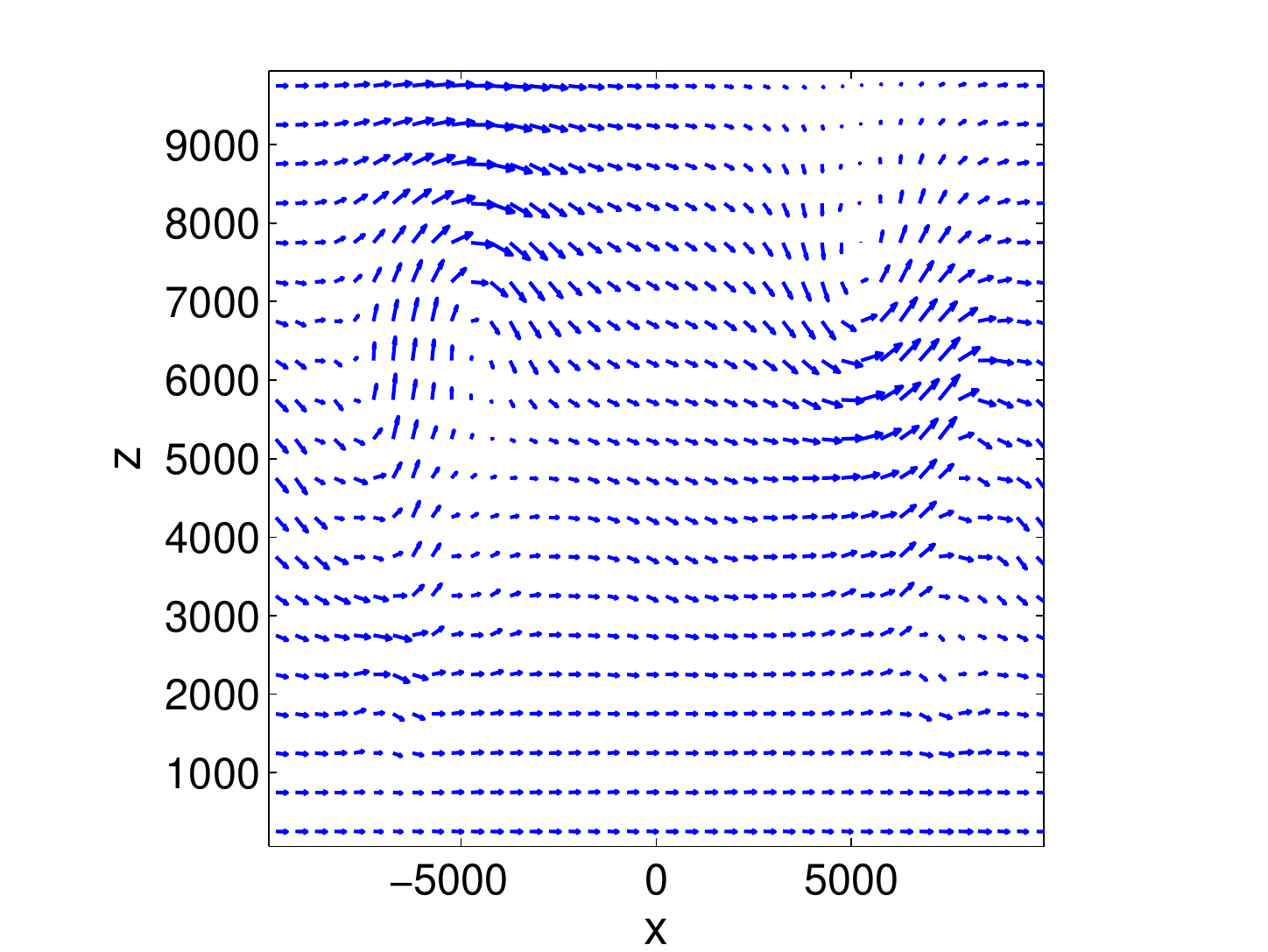}
\end{minipage}
\label{t22}
\end{figure}

\begin{figure}
\caption{Normalized energy (with respect to the initial total value) of the system for the hot and cold bubbles test case.}
\centering
\includegraphics[height=0.3\textheight, width=\textwidth]{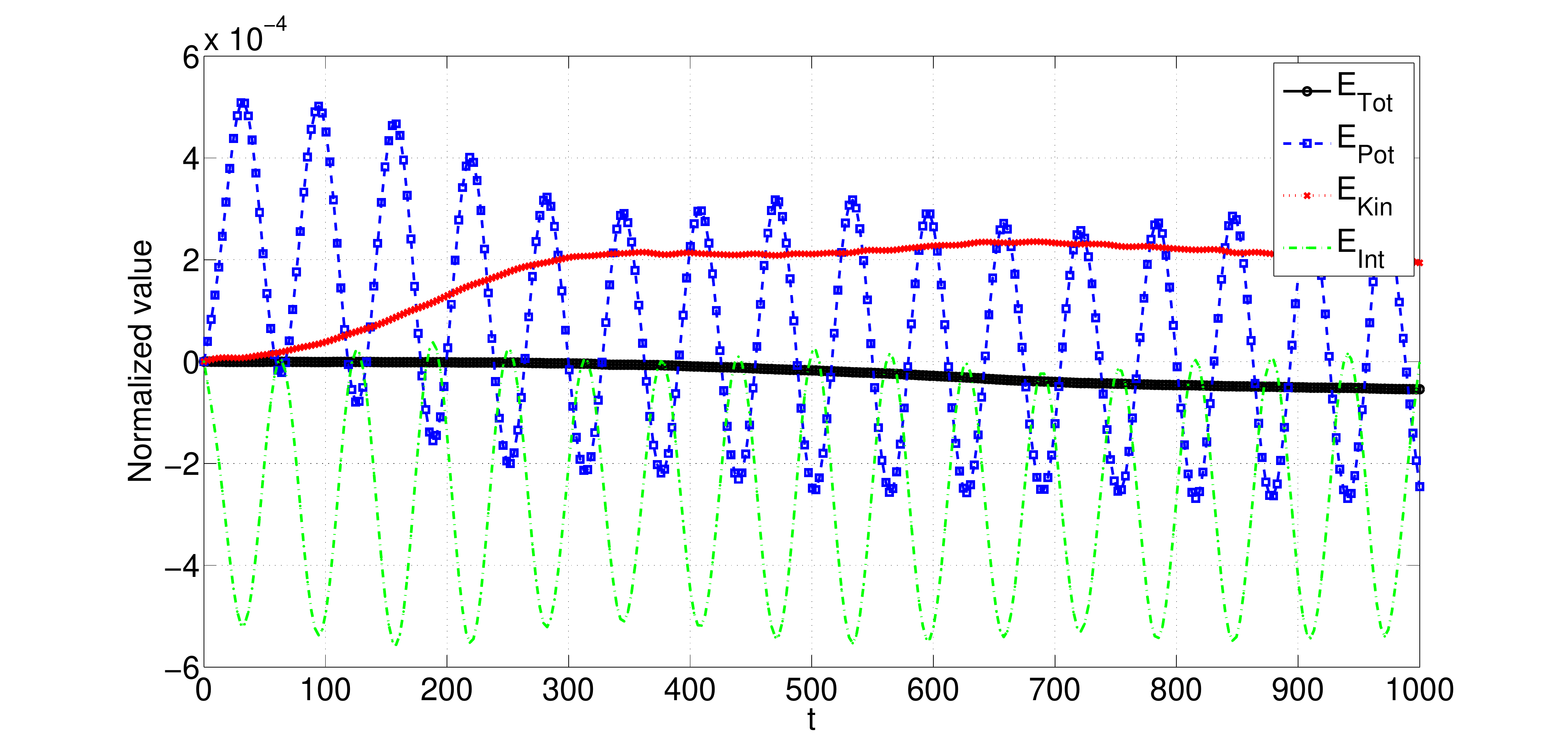}
\label{energyt2}
\end{figure}

We turn our attention to a variation of the previously presented case. Similar tests have being previously performed in \cite{robert}; we have modified it in order to consider the same scales as in our first case, but also to include some kinetic effects to test symmetry. The perturbation in the potential temperature contains a hot but also a cold bubble; the expressions are given by
\begin{eqnarray}
\theta&=&\bar\theta+\theta'_1+\theta'_2,\\
\theta'_1&=&\begin{cases}10\cos\left(\frac{\pi L_1}{2}\right) & L_1\leq 1,\\ 0 & i.o.c.\end{cases},\quad
\theta'_2=\begin{cases}-15\cos\left(\frac{\pi L_2}{2}\right) & L_2\leq 1,\\ 0 & i.o.c.\end{cases},\\
L_1&=&\frac{1}{2000}\sqrt{x^2+(z-2000)^2},\quad
L_2=\frac{1}{2000}\sqrt{x^2+(z-8000)^2}.\\
\end{eqnarray}
We also expect in this case that the hot bubble will rise, but on the other hand the cold bubble should fall, and as they are placed along the same vertical axis, they will collide and interact creating eddy patterns. We consider solid wall boundary conditions at the bottom and the top of the domain, and periodic boundary conditions at the lateral extremes, as we also include an initial horizontal velocity $u=20[ms^{-1}]$ to test the capacity of the scheme to preserve symmetries in the presence of horizontal translation.
Figure \ref{initt2} illustrates the initial setting, figures \ref{t2} and \ref{t22} exhibiting the evolution of the system with a resolution of $\dx=\dz=125[m]$ up to 1000$[s]$. The results reflect the proper physical solution, the rise of the hot bubble together with the fall of the cold bubble until the collision is produced, and the consequent eddy generation due to the interaction of the perturbations. The rotational behavior can be also noticed in the velocity plots, where eddies are also formed. The final state at $1000[s]$ is symmetric with respect to the axis $x=0$, which is coherent with the initial horizontal velocity. No spurious oscillations are observed. The energy plot in figure \ref{energyt2} shows an increment in the kinetic energy of the system, interaction between internal and potential energy, but shows almost no artificial diffusion for the total energy of the system.

\subsection{Density current}
The third case is a popular test case in atmospheric modelling (see \cite{straka,giraldo,carpenter,ahmad}); in the domain $\Omega=[0,\, 20000]\times[0,\,6000]$, with a system initially at rest, a cold bubble perturbation is added to the reference state,
\begin{equation}
\theta'=\begin{cases}-7.5\left(\cos\left(\pi L\right)+1\right) & L\leq 1,\\ 0 & i.o.c.\end{cases},\quad L=\sqrt{\left(\frac{x}{4000}\right)^2+\left(\frac{z-2000}{2000}\right)^2}.
\end{equation}
With solid boundary conditions, the cold bubble drops until it hits the boundary, generating horizontal shear displacement, end eventually developing Kevin-Helmholtz rotors. Final simulation time is $t=900[s]$. It has been previously reported in \cite{straka,giraldo} that most numerical methods require to take into account viscosity effects in order to generate a grid-convergent solution; nevertheless, there exists methods that manage to accurately simulate this test case with a inviscid set of equations \cite{knothmwr}. When a viscous stress term is included in the scheme, a reconstruction step has to be included in the source term approximation and second spatial derivatives are computed from the reconstructed polynomials. In our case, we include simulations for the inviscid model, but we also add a viscous source term $\V$ to the r.h.s. of the system (\ref{euler3d}), with
\begin{equation}
\V=\left[\begin{array}{c}0\\ \rho K (\partial_x^2 u+\partial_z^2 u)\\ \rho K (\partial_x^2 w+\partial_z^2 w)\\ \rho K (\partial_x^2 \theta+\partial_z^2 \theta)\end{array}\right], \quad K=75[m^2s^{-1}].
\end{equation}

\begin{figure}
\caption{Density current test case. Results after 15 min. From top to bottom: potential temperature colormap with $\Delta x=\Delta z=200$, 100 and 50[m]. Left: experiments without viscosity. Right: results with Fickian viscosity with parameter $\nu=75 [m^2\, s^{-1}]$.}
\begin{minipage}[b]{0.5\linewidth}
\centering
\includegraphics[height=0.75\textheight, width=\textwidth]{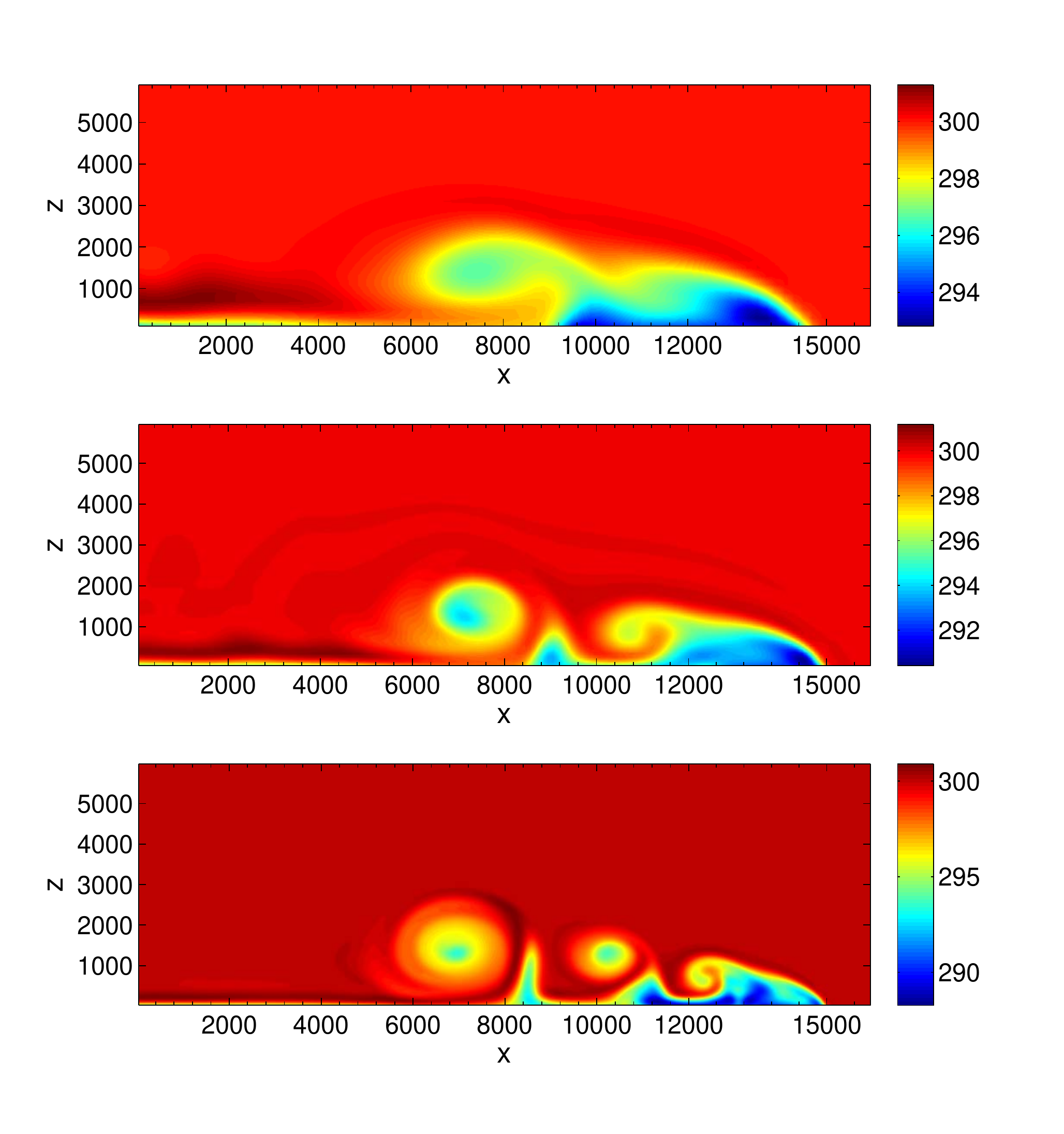}
\end{minipage}
\hspace{0.0cm}
\begin{minipage}[b]{0.5\linewidth}
\centering
\includegraphics[height=0.75\textheight, width=\textwidth]{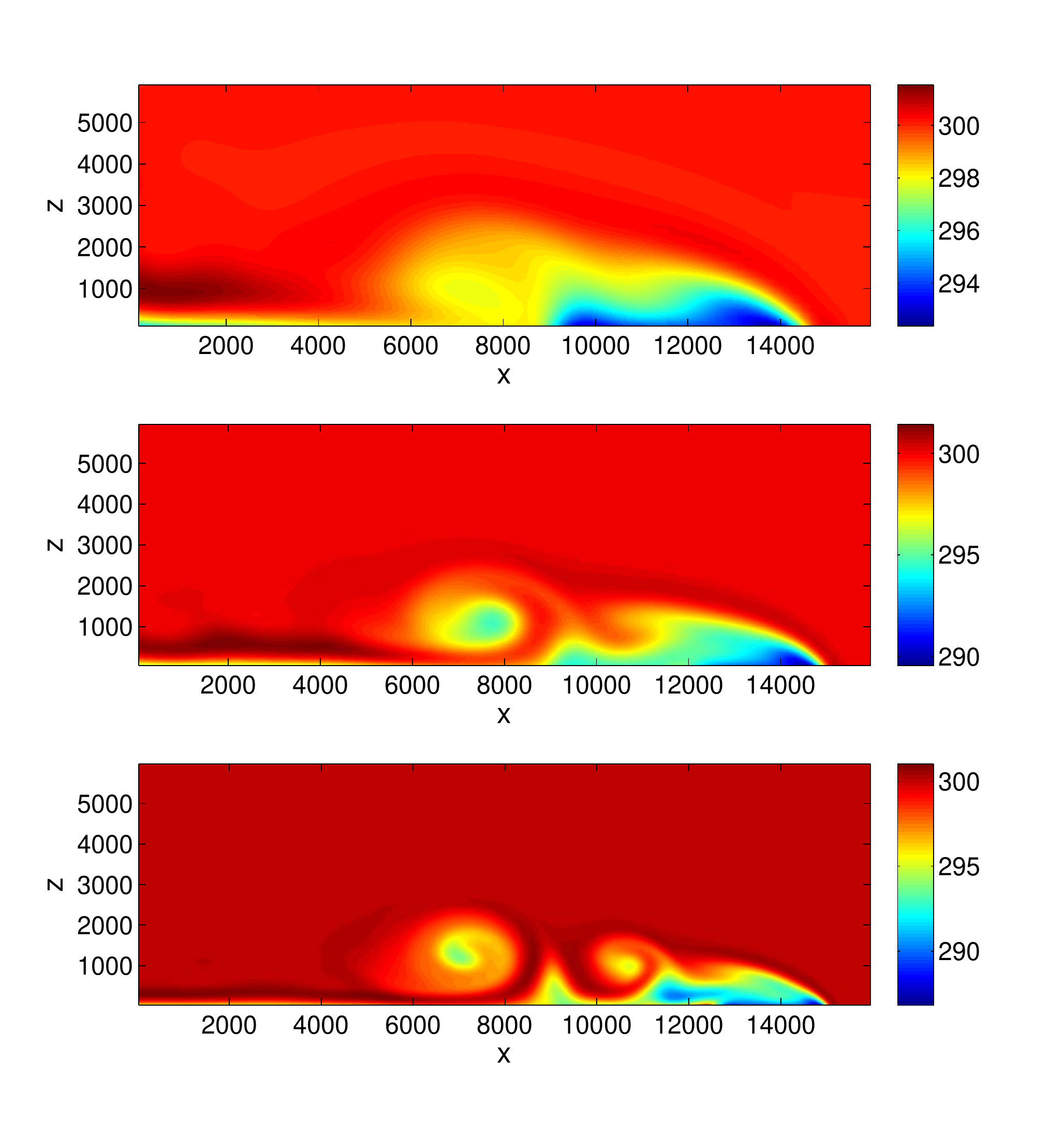}
\end{minipage}
\label{ttest7}
\end{figure}

\begin{figure}
\caption{Density current test case. Results after 15 min and $\Delta x=\Delta z=50$[m], $400\times120$ elements. Top: horizontal velocity. Bottom: vertical velocity. Left: experiments without viscosity. Right: results with Fickian viscosity with parameter $\nu=75 [m^2\,s^{-1}]$.}
\begin{minipage}[b]{0.5\linewidth}
\centering
\includegraphics[height=0.5\textheight, width=\textwidth]{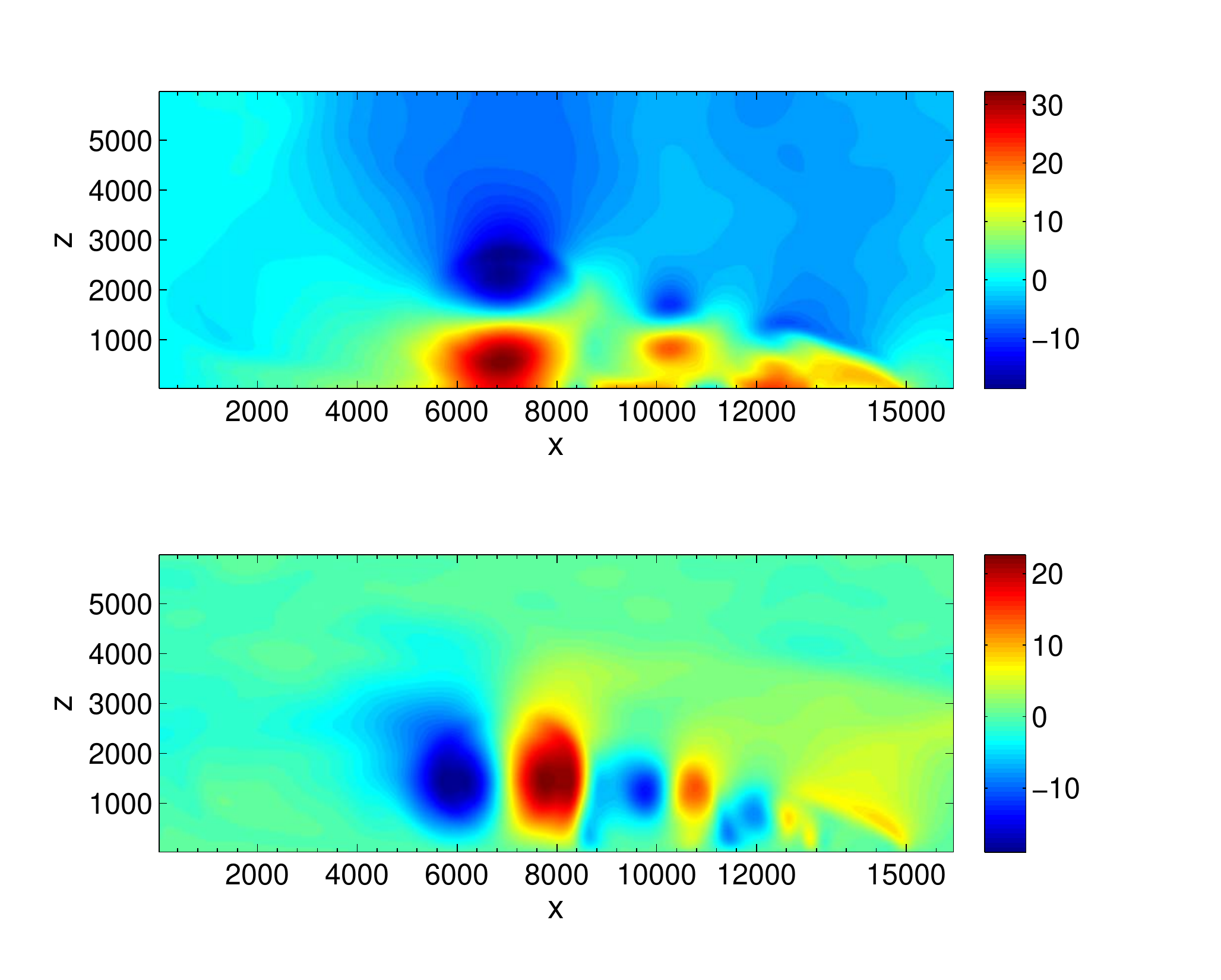}
\end{minipage}
\hspace{0.0cm}
\begin{minipage}[b]{0.5\linewidth}
\centering
\includegraphics[height=0.5\textheight, width=\textwidth]{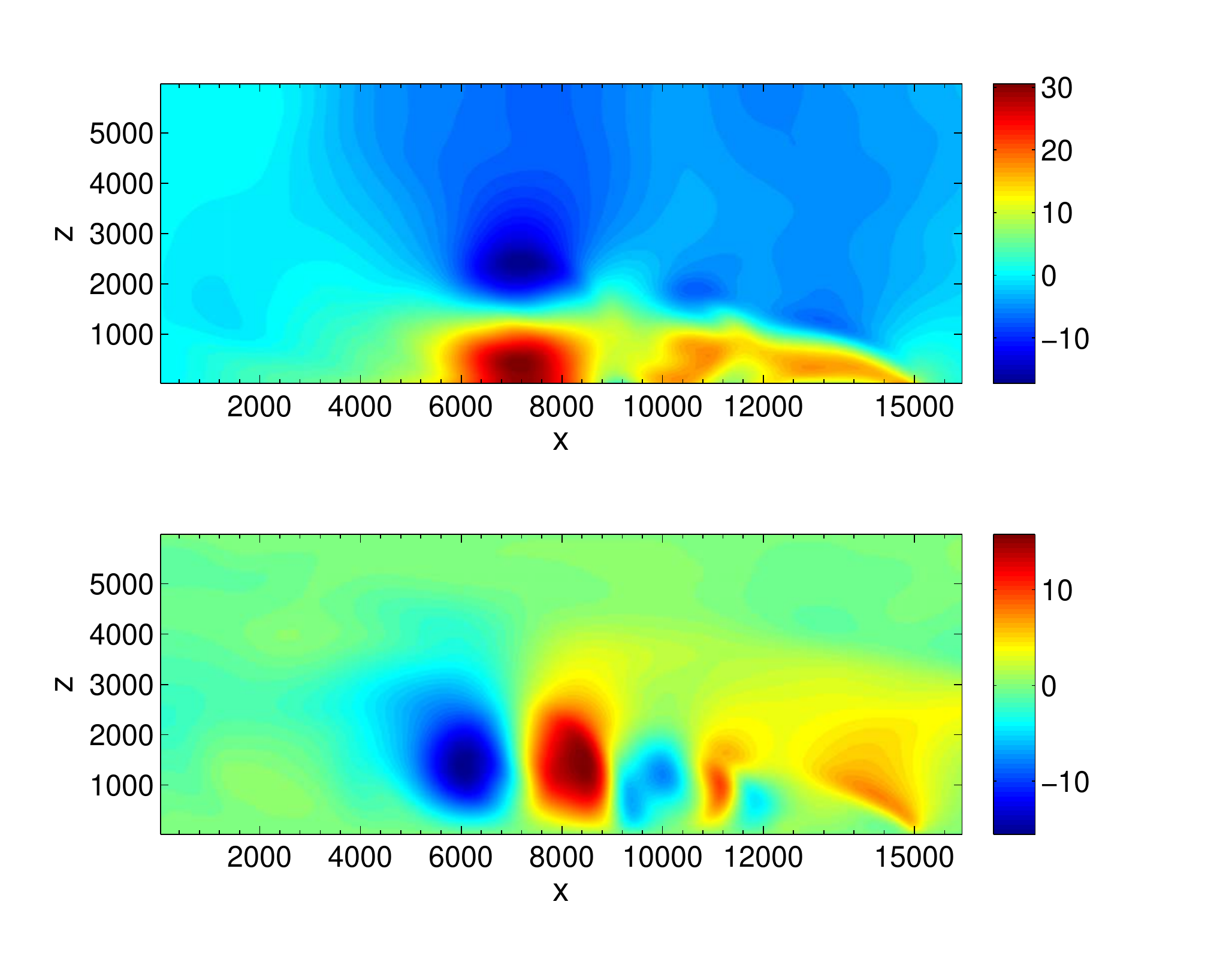}
\end{minipage}
\label{vetest7}
\end{figure}

\begin{figure}
\caption{Normalized energy (with respect to the initial total value) of the system for the density current test case without viscosity.}

\centering
\includegraphics[height=0.3\textheight, width=\textwidth]{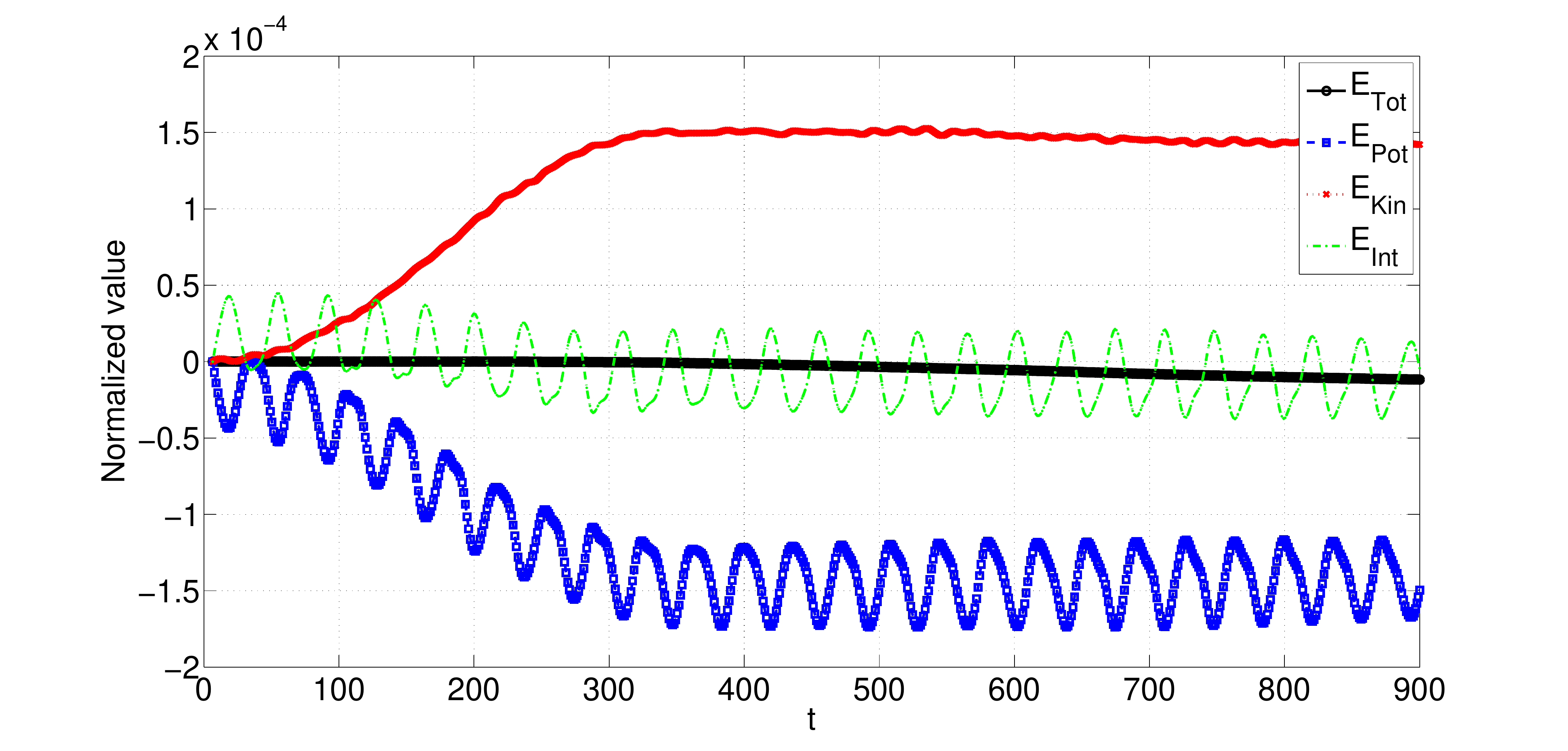}
\label{etest7}
\end{figure}

\begin{figure}
\caption{Normalized energy (with respect to the initial total value) of the system for the density current test case with viscosity.}

\centering
\includegraphics[height=0.3\textheight, width=\textwidth]{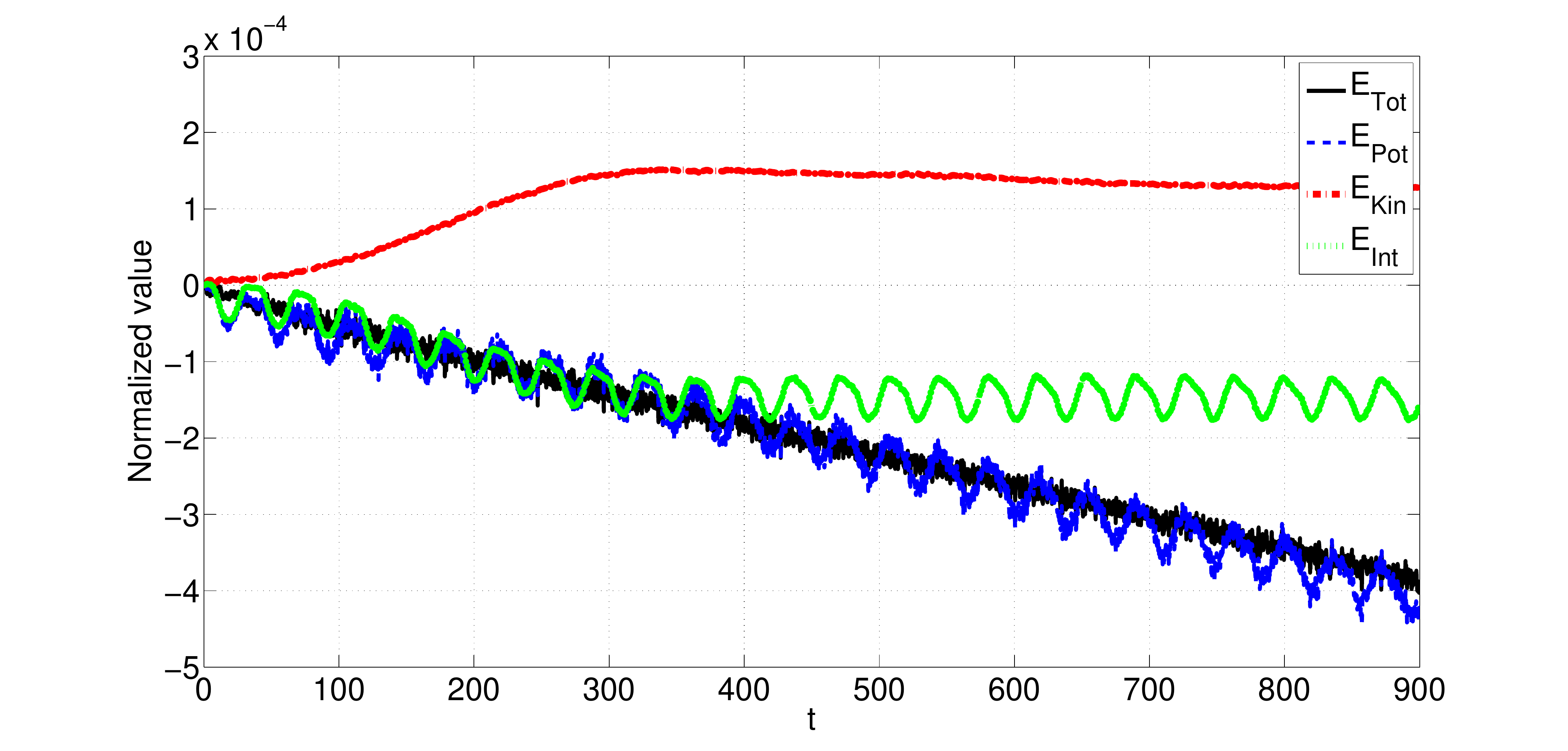}
\label{evtest7}
\end{figure}

We first study the effect of resolution in both inviscid and viscous simulations. In figures \ref{ttest7} and \ref{vetest7} it can be seen final time results for variable resolution. In both cases, there is a clear convergent behavior of the solution, and increasing resolution provides a better insight into the development of Kevin-Helmhotlz rotors; at the highest tested resolution, with $\dx=\dz=50[m]$, it can be clearly appreciated the generation of three rotors, which is similar to the results obtained in the aforementioned references for this test case for viscous simulations (including the diffused aspect of the eddy nearest to the front). Table \ref{tab:tdc} shows the extreme values obtained for both simulations at high-resolution, which are in accordance to the range of values previously obtained by other authors; in particular, the front location, which is a relevant quantity in this test case, coincides the results published in \cite{giraldo,straka,ahmad}. The energy plots shown in figure \ref{etest7} and \ref{evtest7} are qualitatively comparable to the one presented in \cite{carpenter} (even though is not exactly the same density current test case), exhibiting a sustained increment in the kinetic energy of the system, while both internal and potential energy decrease throughout the simulation. There are some differences though, when it comes to analyze the conservation of the total energy: the inviscid system is completely conservative, while the inclusion of diffusion alters this property and decreases in time.

\begin{table}[H]
\centering
\caption{Extreme values for the density current test case with $\dx=\dz=50[m]$ at $t=900[s]$.}\label{tab:tdc}\vskip 3mm
\begin{tabular}{c c c c c c c c}
\hline\\
Model& $\theta'_{min}$  & $\theta'_{max}$ & $u_{min}$ & $u_{max}$ & $w_{min}$ & $w_{max}$ & Front location   \\ [0.5ex]
\hline\\
Inviscid & -11.7296 & 0.8950 & -18.6235 &32.2655 &-18.9703 & 22.6452 & 1.498$\times 10^4$ \\
Viscous&  -13.1851& 1.0209 & -17.2731& 30.5931 & -15.5440 & 15.7718  & 1.503$\times 10^4$\\
\hline
\end{tabular}
\end{table}

\subsection{Convective bubble in a stable atmosphere}
\begin{figure}
\caption{Convective bubble in a stable atmosphere. Results after 10 min. From top to bottom: potential temperature, horizontal velocity and vertical velocity. Left: low resolution results with $\Delta x=\Delta z=500$[m]. Right: high resolution results with $\Delta x=\Delta z=250$ [m]. }
\begin{minipage}[b]{0.5\linewidth}
\centering
\includegraphics[height=0.75\textheight, width=\textwidth]{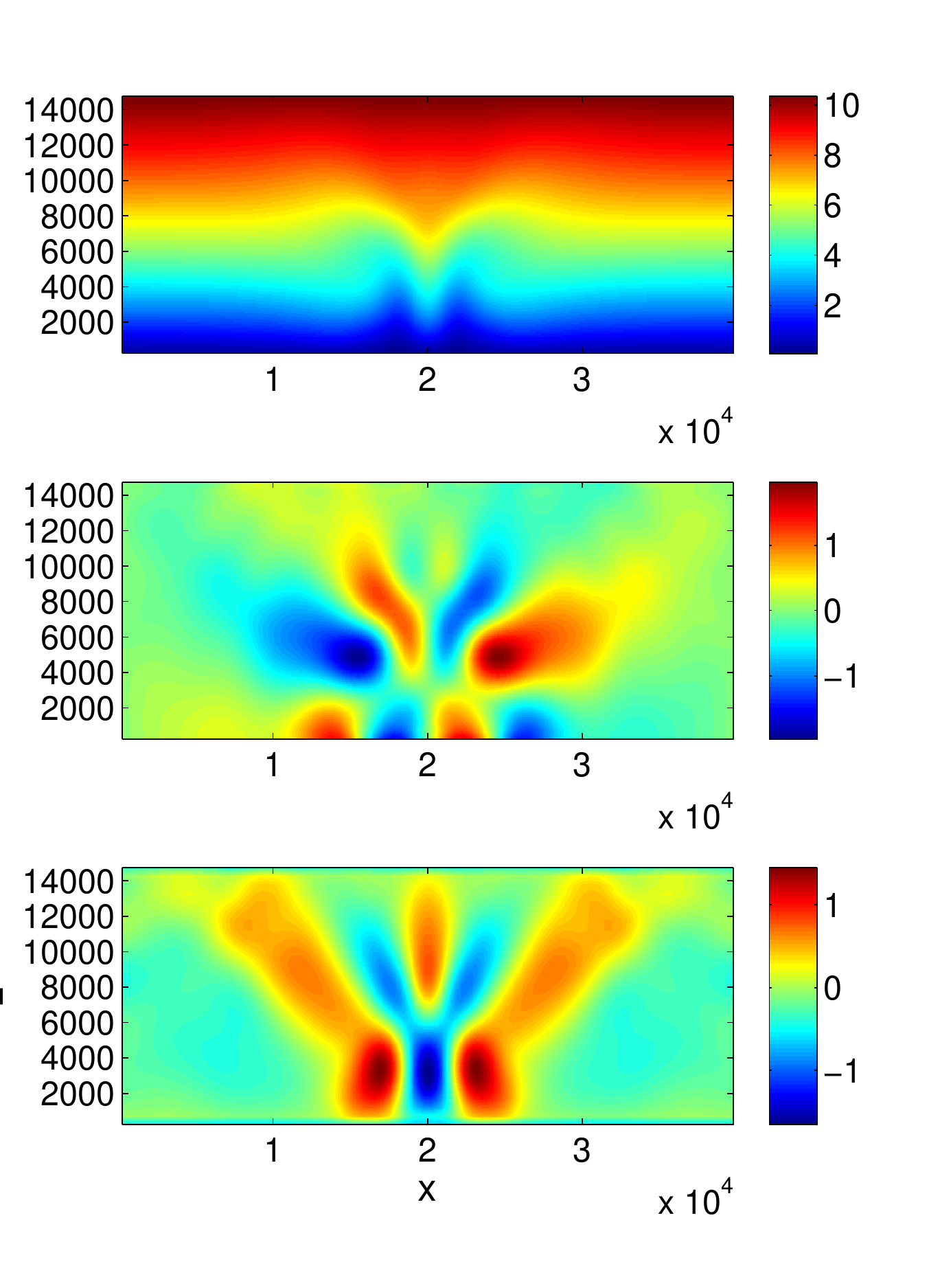}
\end{minipage}
\hspace{0.0cm}
\begin{minipage}[b]{0.5\linewidth}
\centering
\includegraphics[height=0.75\textheight, width=\textwidth]{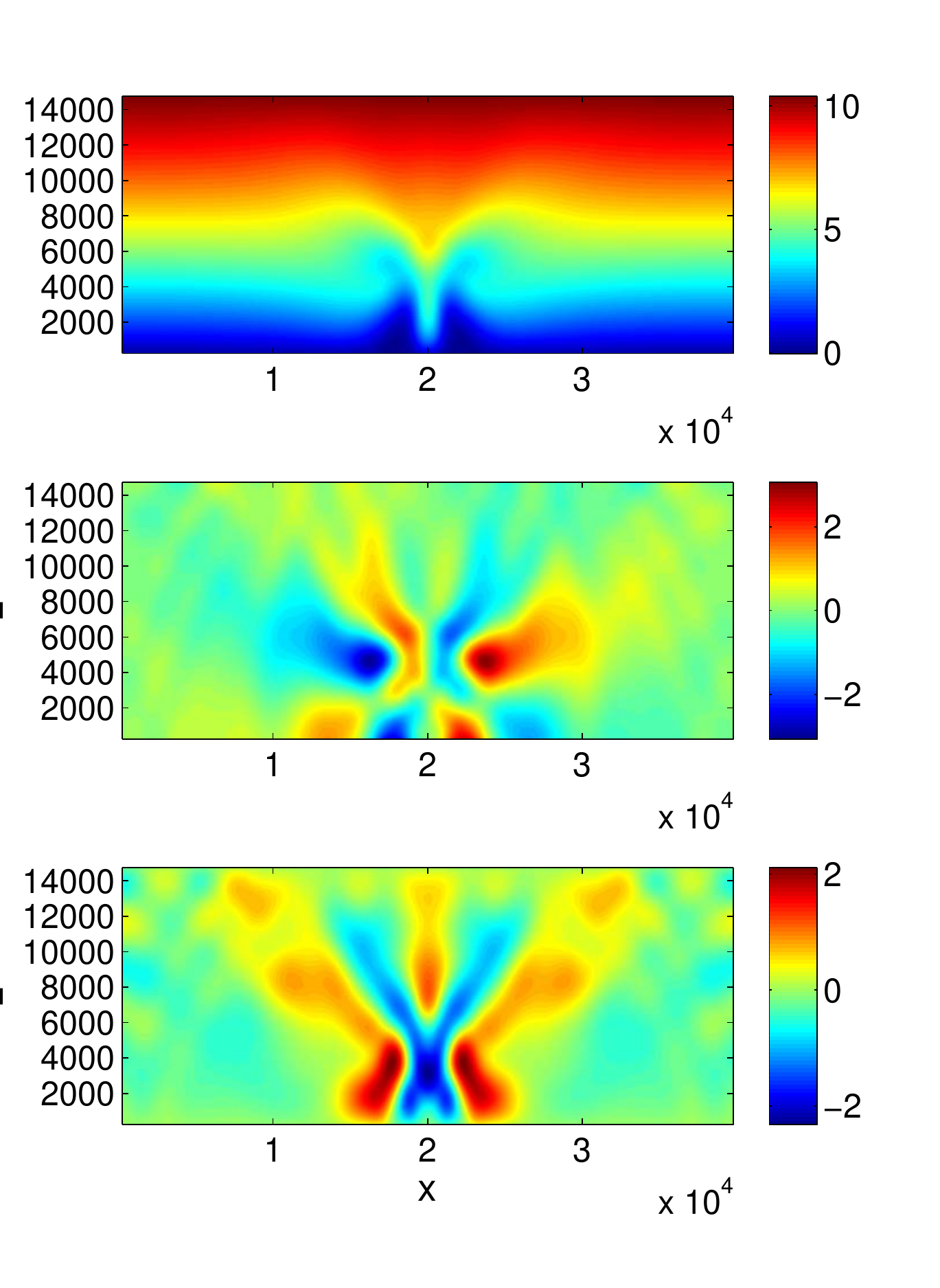}
\end{minipage}
\label{fig:t6}
\end{figure}

\begin{figure}
\caption{Normalized energy (with respect to the initial total value) of the system for convective bubble in a stable atmopshere.}
\centering
\includegraphics[height=0.3\textheight, width=\textwidth]{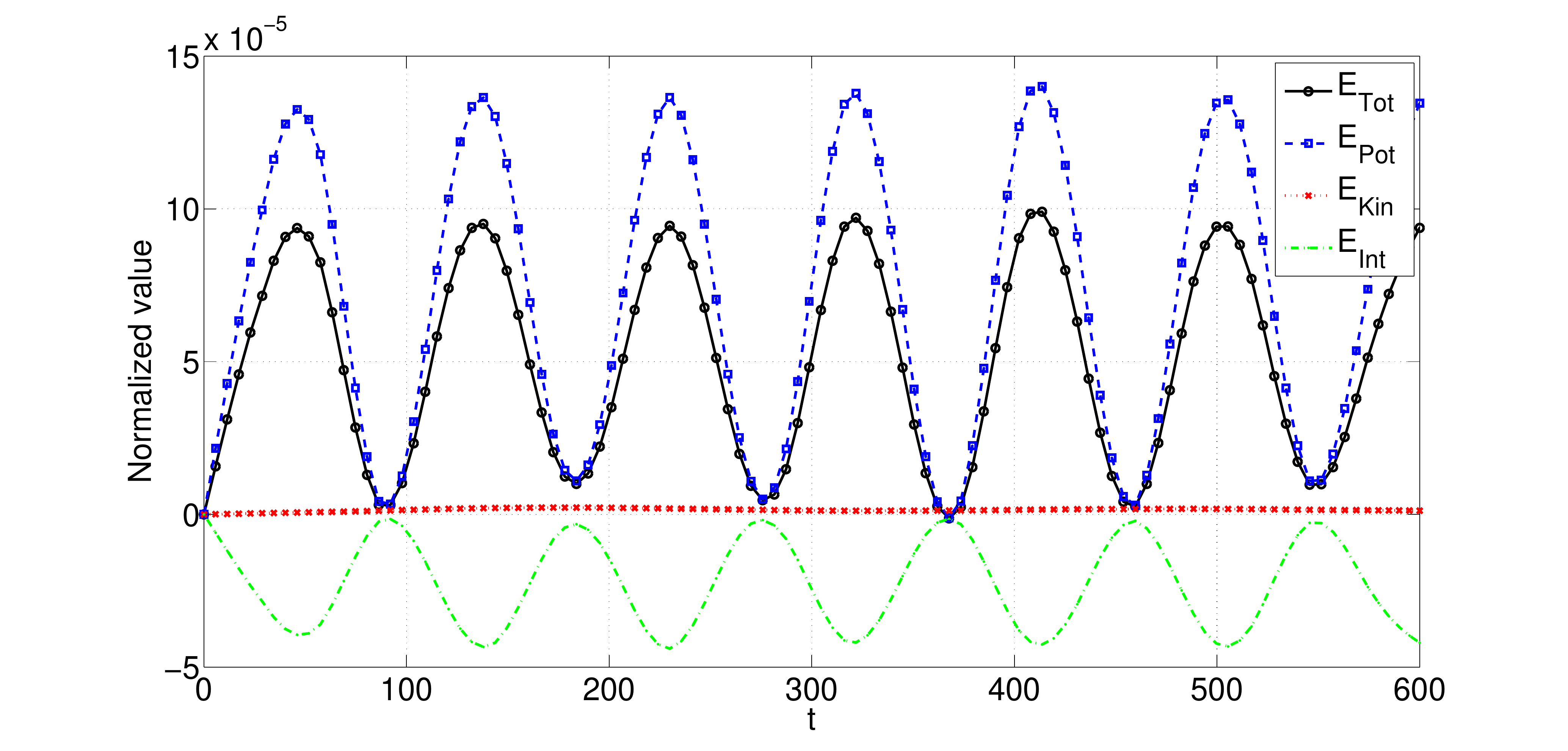}
\label{fig:et6}
\end{figure}

We conclude our study with a test that has been previously presented in \cite{mendez}. In the domain $\Omega=[0,\,40000]\times[0,\,15000]$, we consider a stable atmosphere with an initial potential temperature vertical gradient of $4[K\,km^{-1}]$ with a mean ground level value of $300[K]$. After vertical integration of the gradient, a positively stratified reference state for the potential temperature is obtained, and density is initialized via hydrostatic balance. We add a warm potential temperature perturbation bubble of the form
\begin{equation}
\theta'=\begin{cases}6.6\cos^2\left(\frac{\pi L}{2}\right) & L\leq 1,\\ 0 & i.o.c.\end{cases},\quad L=\frac{1}{2500}\sqrt{x^2+(z-2750)^2}.
\end{equation}
Simulation time is set to $t=600[s]$, we use solid wall boundary conditions at the top and the bottom of the domain, and open boundary conditions for the lateral extremes. We execute model runs with $\dx=\dz=500 [m]$ and $\dx=\dz=250 [m]$.
Figure \ref{fig:t6} shows high and low-resolution results at final time of simulation. Low-resolution results are in accordance to what is presented in \cite{mendez}; since the temperature increases with altitude, a mitigation of the buoyancy is expected while horizontal spreading of the perturbation occurs. High-resolution experiments exhibits the same behavior, although there is a variation in the extremal values as it can be seen in table \ref{tab:t6}. In both cases symmetry with respect to $x=0$ is preserved. Regarding energy conservation, figure \ref{fig:et6} shows that total energy is preserved while potential, internal and kinetic energy oscillate with respect the initial state.

\begin{table}[H]
\centering
\caption{Extreme values for the convective bubble in a stable atmosphere test case at $t=600[s]$.}\label{tab:t6}\vskip 3mm
\begin{tabular}{c c c c c c c}
\hline\\
Resolution& $\theta'_{min}$  & $\theta'_{max}$ & $u_{min}$ & $u_{max}$ & $w_{min}$ & $w_{max}$   \\ [0.5ex]
\hline\\
$\dx=500$ [m]& 0.0748 & 10.3438 & -1.9717 & 1.9717 & -1.6539 & 1.4501 \\
$\dx=250$ [m]& -0.0251 & 10.3890 & -3.0724 & 3.0685 & -2.3212 & 2.1149 \\
\hline
\end{tabular}
\end{table}

\section{Summary and outlook}
We have presented a second-order, non-oscillatory scheme for the resolution of relevant advective and convective atmospheric phenomena. The method makes use of a WENO reconstruction procedure for accurate extrapolation of boundary and inner cell values, together with a centered-limited approach for the flux calculation. Even though, for this class of problems, schemes based on upwinding considerations are usually preferred over centered approaches, the proposed scheme performs well in an extensive set of test problems. The theoretically expected second order was reached whenever an analytic solution was available, and in its absence the scheme proved to be grid-convergent with a performance similar to the currently available algorithms. The accuracy of the method can be numerically increased by switching the flux parameter $\omega$, with a cost associated to the time stepping in order to preserve monotonicity. The scheme also showed robustness with respect to the limiter choice in a sharp front, which also illustrates the robustness of the reconstruction procedure. In the convective experiments, the scheme managed to reproduce the correct physical behavior, while tracking fronts in the correct position and without generating spurious oscillations. It has also showed robustness with respect of the inclusion of viscosity: unless it is physically relevant, the current evidence seems to indicate that the method is able to perform in a consistent way without the need of viscosity, which is a consequence of the continuous enforcement of the non-oscillatory character of the scheme, present in both time and spatial discretizations.

We point some open issues to address. The removal of the splitting approach, by trying to incorporate source terms effects into the flux calculation could make the method far more efficient. In the splitting, context, other alternatives for viscosity calculation should be explored, as in the current scheme additional reconstruction steps (one of the most expensive parts of the code) are required. Nevertheless, we conclude that the studied scheme possess a great level of applicability in the atmospheric modelling framework; further enhancements on its formulation and performance could increase its potential as a robust, high-resolution method.





\bibliographystyle{model3-num-names}
\bibliography{bibtesis}



 %

%
%
\end{document}